\documentclass[12pt]{article}
\usepackage[latin1]{inputenc}
\usepackage{amsfonts,amssymb,amsmath, epsfig}
\usepackage{color,graphicx,graphics}
\usepackage{amsthm}
\usepackage{amsmath,amstext,amssymb,amsfonts, amscd}
\usepackage[lofdepth,lotdepth]{subfig}
\usepackage{xcolor}
\usepackage{hyperref}          

\def\Box{\leavevmode\vbox{\hrule
     \hbox{\vrule\kern4pt\vbox{\kern4pt}%
           \vrule}\hrule}}
\def\blackbox{\leavevmode\vrule height 5pt width 4pt depth 0pt\relax}
\def\endproof{\null\hfill {$\blackbox$}\bigskip}

\def\paragraph#1{{\bf #1\ }}

\newtheorem{lemma}{Lemma}[section]  

\newtheorem{theorem}[lemma]{Theorem}

\newtheorem{proposition}[lemma]{Proposition}

\newtheorem{remark}{Remark}[section]

\title{Radial Laplacian on rotation groups} 
\author{Pierre Degond} 
\date{} 
\begin{document}

\maketitle

$\mbox{}$

\vspace{-1 cm}

\begin{center}
Institut de Math\'ematiques de Toulouse ; UMR5219 \\
Universit\'e de Toulouse ; CNRS \\
UPS, F-31062 Toulouse Cedex 9, France \\
pierre.degond@math.univ-toulouse.fr
\end{center}

\vspace{0. cm}

\begin{center}
{\em In fond memory of Sao Carvalho, estimated colleague and friend}
\end{center}

\vspace{0. cm}
\begin{abstract}
The Laplacian on the rotation group is invariant by conjugation. Hence, it maps class functions to class functions. A maximal torus consists of block diagonal matrices whose blocks are planar rotations. Class functions are determined by their values of this maximal torus. Hence, the Laplacian induces a second order operator on the maximal torus called the radial Laplacian. In this paper, we derive the expression of the radial Laplacian. Then, we use it to find the eigenvalues of the Laplacian, using that characters are class functions whose expressions are given by the Weyl character formula. Although this material is familiar to Lie-group experts, we gather it here in a synthetic and accessible way which may be useful to non experts who need to work with these concepts. 
\end{abstract}

\medskip
\noindent
{\bf Acknowledgements:} PD holds a visiting professor association with the Department of Mathematics, Imperial College London, UK.

\medskip
\noindent
{\bf Key words: } Representations, characters, eigenfunctions

\medskip
\noindent
{\bf AMS Subject classification: } 22E30, 22E46, 22E70, 81R05
\vskip 0.4cm

\setcounter{equation}{0}
\section{Introduction}
\label{intro}

Rotation groups are of fundamental physical importance because they encode symmetries of the underpinning space. Recently, they have emerged as important players in the description of large systems of collectively moving self-propelled particles. Such systems are frequently encountered  in life science, such as swarming bacteria, migrating cells or schooling fish. While motion of physical particles is determined by their momenta which span a flat space, collectively moving particles offer a greater variety of situations. It often happens that particle dynamics is described in terms of elements of a Lie group or of a homogeneous space (i.e. a quotient of a Lie group by one of its subgroups). 

One of the earliest examples of such models is the Vicsek model \cite{vicsek1995novel}, where the particles are moving with constant speed. Hence, their dynamics is determined by the direction of their velocity which belongs to the homogeneous space ${\mathbb S}^{n-1} \cong \textrm{SO}_n{\mathbb R}/\textrm{SO}_{n-1}{\mathbb R}$. More recently, a model in which particle dynamics is determined by an element of the rotation group $\textrm{SO}_{n}{\mathbb R}$ has been proposed \cite{degond2017new}. Several variants of this model have been explored \cite{degond2021body, Degond_eal_JNLS20, Degond_etal_proc19, Degond_etal_MMS18}. One may find other collective dynamics models involving the rotation group in \cite{fetecau2021emergent}, or involving other groups or homogeneous spaces, such as ${\mathbb S}^1 \times {\mathbb R}$ in \cite{degond2011macroscopic},  ${\mathbb S}^n \times {\mathbb S}^1$ in~\cite{degond2022topological} or the unitary group $\textrm{SU}_n$ in \cite{golse2019mean, ha2017emergent}. More generally, the interplay between collective dynamics and geometry is an active subject, see e.g. \cite{ahn2021emergent, fetecau2019self}. 

Here, we focus on the rotation group in relation to the model exposed in \cite{degond2017new}. In this paper, a particle system for rigid bodies interacting through body attitude alignment was proposed and a kinetic model was derived when the number of particles is large. Finally, a hydrodynamic model was deduced from the kinetic model assuming that the alignment interaction is strong. At the kinetic level, the interaction operator is of Fokker-Planck type and thus, involves the Laplacian on $\textrm{SO}_n{\mathbb R}$. The derivation of the hydrodynamic model made in \cite{degond2017new} strongly used the peculiar structure of $\textrm{SO}_3{\mathbb R}$ (Rodrigues formula or the quaternion representation like in \cite{Degond_etal_MMS18}) and thus, was restricted to the case of dimension $n=3$. A further step towards an arbitrary dimension $n$ was made in \cite{degond2021body}, but in this case, the interaction operator was replaced by a relaxation operator of BGK type, thus avoiding to deal with the Laplacian on $\textrm{SO}_n{\mathbb R}$. The goal of this paper is to collect the appropriate mathematical framework in view of extending the result of \cite{degond2017new} to arbitrary dimensions. Here we focus on the rotation group. Its use in the model of \cite{degond2017new} will be the subject of future work. 

Although the material presented here is well-known to Lie-group experts, we were not able to find it anywhere exposed in an accessible fashion to non experts. In particular, we could not find anywhere the expression of the radial Laplacian as exposed in Section \ref{sec_rad_lap} below. The present work owes a lot to the presentation made in \cite{faraud2008Analysis} but \cite{faraud2008Analysis} restricts itself to the unitary group $\textrm{SU}_n$ so we had to reprocess all the details for the rotation groups. By applying the radial Laplacian to the characters of the irreducible representations (which are known thanks to the Weyl character formula), one determines all the eigenvalues of the Laplacian, as explained in Section \ref{sec_prop}. Ref. \cite{fulton2013representation} gives a comprehensive introduction to Lie group and Lie algebra representations. Another point of view, focused on compact Lie group, can be found in \cite{Simon}. Finally, the role of the radial Laplacian in the spectral theory of the Laplacian is treated in \cite{gurarie2007symmetries, helgason2022groups}, however in an abstract way which makes it hardly tractable in practice. Of course, Lie group theory is a very mature subject and there are numerous books and treatises on it. The above is just a small selection. Let us also mention the classical books \cite{chevalley2018theory, weyl1946classical}. Other useful references are (without being exhaustive) \cite{adams1982lectures, brocker2013representations, bump2004lie, duistermaat2012lie, hall2013lie, humphreys2012introduction, knapp1996lie, sepanski2007compact}. 

The organization of this paper is as follows. Section \ref{sec_lap} introduces the basic notations and facts about the rotations groups, the Laplacian operator acting on them and the radial Laplacian. Section \ref{sec_rad_lap} is devoted to the computation of the radial Laplacian while Section \ref{sec_prop} develops the spectral theory of the Laplacian and radial Laplacian. Finally, a conclusion is drawn in Section \ref{sec_conclu}. Two appendices provide proofs of auxilliary results.

\setcounter{equation}{0}
\section{Laplacian on rotation groups}
\label{sec_lap}

The following presentation is succinct. We refer to e.g. \cite{faraud2008Analysis} for detail. We denote by $\textrm{SO}_n{\mathbb R}$, the group of rotations of ${\mathbb R}^n$, i.e. the subset of the space $M_n{\mathbb R}$ of $n \times n$ real matrices defined by 
$$ \textrm{SO}_n{\mathbb R} = \{ A \in M_n{\mathbb R} \, \, | \, \, A A^T = A^T A = \textrm{I}, \quad \det A = 1 \}, $$
where $A^T$ denotes the transpose of $A$, $\textrm{I}$ is the identity matrix of ${\mathbb R}^n$ and $\det A$ is the determinant of $A$. The space $M_n{\mathbb R} \approx {\mathbb R}^{2n}$ is a vector space which can be endowed with a Euclidean structure by means of the following inner product 
\begin{equation}
A \cdot B = \frac{1}{2} \textrm{Tr} \{ A^T B \}, \quad \forall A, \, B \in M_n{\mathbb R}, 
\label{eq:inner_product}
\end{equation}
where $\textrm{Tr} A$ denotes the trace of the matrix $A$. We note that \eqref{eq:inner_product} is the usual Frobenius inner product multiplied by the factor $\frac{1}{2}$. This turns out to be the convenient definition of a matrix inner product when dealing with rotations groups. The set $\textrm{SO}_n{\mathbb R}$ is at the same time a group for matrix multiplication and an embedded manifold into $M_n{\mathbb R}$ where matrix multiplication is $C^\infty$. It is therefore a Lie group. Its tangent space at $\textrm{I}$ is the Lie algebra $\mathfrak{so}_n{\mathbb R}$ of antisymmetric matrices endowed with the Lie bracket $[P,Q] = PQ-QP$, for all $P, \, Q \in 
\mathfrak{so}_n{\mathbb R}$. Its tangent space at any $A \in \textrm{SO}_n{\mathbb R}$ is given by $T_A \textrm{SO}_n{\mathbb R}= A \, \mathfrak{so}_n{\mathbb R} = \mathfrak{so}_n{\mathbb R} \, A$. From now on, we will abbreviate $T_A \textrm{SO}_n{\mathbb R}$ into simply $T_A$. As a Lie group, $\textrm{SO}_n{\mathbb R}$, can be endowed with a unique (up to a multiplicative factor) Haar measure, i.e. a measure which is invariant by left multiplication. Being compact, $\textrm{SO}_n{\mathbb R}$ is a modular group. Thus, the Haar measure is also right invariant, as well as invariant by group inversion and hence, by transposition. Also, the Haar measure can be normalized, and the resulting normalized Haar measure is unique. 

As an embedded manifold in $M_n{\mathbb R}$, $\textrm{SO}_n{\mathbb R}$ is also endowed with a Riemannian inner-product and metric (for Riemannian geometry, we refer to \cite{do1992riemannian, gallot1990riemannian}). For two elements $AP$, and $AQ$ in $T_A$, with $P$ and $Q$ in $\mathfrak{so}_n{\mathbb R}$, the Riemannian inner-product is given by $(AP, AQ) = AP \cdot AQ = P \cdot Q$. On oriented manifolds, the Riemannian metric gives rise to a Riemannian volume form $\omega$ and thus to a positive measure written $d \omega$ by abuse of notation. As any Lie group is orientable \cite[Section 4.3, Example (a)]{warner1983foundations}, $\textrm{SO}_n{\mathbb R}$ is thus endowed with a Riemannian volume form~$\omega$. On general Lie groups, the Riemannian measure $d \omega$ and the Haar measure may not be multiple of each other. However, they are indeed multiple of each other on $\textrm{SO}_n{\mathbb R}$ as stated in the following Lemma, whose proof is shown in Appendix \ref{secapp_proportional}: 

\begin{lemma}
Let $d\omega$ and $d \mu$ be the Riemannian measure and the Haar measure on $\textrm{SO}_n{\mathbb R}$ respectively. Then, there exists a constant $C_n >0$ such that  $d \omega = C_n d\mu$. 
\label{lem:vol_form}
\end{lemma}

On a Riemannian manifold, we can define the gradient and divergence operators. The gradient of a smooth map $\textrm{SO}_n{\mathbb R} \to {\mathbb R}$, $A \mapsto f(A)$ at $A$ is an element $\nabla f (A)$ of $T_A$ defined by 
\begin{equation}  
\nabla f (A) \cdot X = df_A (X), \quad \forall X \in T_A, 
\label{eq:gradient}
\end{equation}
where $df_A$ is the differential of $f$ at $A$. We recall that $df_A$ is a linear map $T_A \to {\mathbb R}$ and $df_A(X)$ is the image of the tangent vector $X$ by this map. A smooth tangent vector field $\phi$ is a smooth map $\textrm{SO}_n{\mathbb R} \to M_n{\mathbb R}$, $A \mapsto \phi(A)$ such that $\phi (A) \in T_A$, for all $A \in \textrm{SO}_n{\mathbb R}$. For instance, the gradient of a smooth map $\textrm{SO}_n{\mathbb R} \to {\mathbb R}$ is a smooth vector field. The divergence $\nabla \cdot \phi$ of the tangent vector field $\phi$ is defined by duality as follows: 
$$ \int_{\textrm{SO}_n{\mathbb R}} \nabla \cdot \phi \, \psi \, d\omega = - \int_{\textrm{SO}_n{\mathbb R}} \phi \cdot \nabla \psi \, d\omega, $$
where $\psi$ is an arbitrary smooth map from $\textrm{SO}_n{\mathbb R}$ to ${\mathbb R}$. Here, we note the use of the Riemannian measure in the integrals. However, because of Lemma \ref{lem:vol_form}, the Riemannian measure $d \omega$ can be replaced by the Haar measure $d \mu$. 

The Riemannian Laplacian $\Delta_M$ is defined as the composition of the gradient and the divergence, i.e. 
$$\Delta_M f = \nabla \cdot ( \nabla f). $$ 
If $\phi$ is a tangent vector field, and $f$ is a real function on $\textrm{SO}_n{\mathbb R}$, it is customary in differential geometry \cite[Section 1.22, formula (6)]{warner1983foundations} to define the map $\phi(f)$,  $\textrm{SO}_n{\mathbb R} \to {\mathbb R}$, by
\begin{equation} 
\phi(f) (A)= : df_A \big(\phi (A) \big) = \nabla f (A) \cdot \phi(A). 
\label{eq:dir_deriv}
\end{equation}
In other words, $\phi(f)(A)$ is the directional derivative of $f$ in the direction of $\phi$ at $A$. For $i$, $j$ in \{1, \ldots, n\}, let $E_{ij} = e_i \otimes e_j \in M_n{\mathbb R}$ where $e_i$ is the $i$-th canonical basis vector of ${\mathbb R}^n$ and $\otimes$ is the tensor product. In other words, $(E_{ij})_{k \ell} = \delta_{ij} \delta_{k \ell}$, where $\delta$ is the Kronecker delta. Then, $\big(\sqrt{2} E_{ij}\big)_{(i,j) \in \{1, \ldots, n\}^2}$ is an orthonormal basis of $M_n{\mathbb R}$ for the inner product~\eqref{eq:inner_product}. Denote by $P_{T_A}$ the orthogonal projection of $M_n{\mathbb R}$ on $T_A$. We can easily check that 
\begin{equation}
P_{T_A} M = A \frac{A^T M - M^T A}{2}, \quad \forall M \in M_n{\mathbb R}, \quad \forall A \in \textrm{SO}_n{\mathbb R}. 
\label{eq:PTA_form}
\end{equation}
Then, by \cite[Theorem 3.1.4]{hsu2002stochastic}, we have 
\begin{equation}
\Delta_M f(A) = 2 \sum_{(i,j) \in \{1, \ldots, n\}^2} (P_{T_A} E_{ij})^2 (f)(A), \quad \forall A \in \textrm{SO}_n{\mathbb R}, 
\label{eq:lapriem_def}
\end{equation}
where the square means that the operator $P_{T_A} E_{ij}$ (see \eqref{eq:dir_deriv} for its definition) is acted twice on $f$. The right-hand side of \eqref{eq:lapriem_def} can be shown to be independent of the choice of the orthonormal basis of $M_n{\mathbb R}$. As is clear from above, this definition of the Laplacian does not depend on the Lie-group structure of $\textrm{SO}_n{\mathbb R}$. 

On Lie groups, there is an alternate concept of the Laplacian (see \cite[Section 8.2]{faraud2008Analysis}). We first remind that the matrix exponential maps $\mathfrak{so}_n{\mathbb R}$ into $\textrm{SO}_n{\mathbb R}$ (and actually onto, but this will not be used). For any $X \in \mathfrak{so}_n{\mathbb R}$, $\rho(X)$ denotes the linear map: $C^\infty(\textrm{SO}_n{\mathbb R}) \to C^\infty(\textrm{SO}_n{\mathbb R})$, $f \mapsto \rho(X) f$ such that 
$$ \big(\rho(X) f \big) (A) = \frac{d}{dt} \big( f (A e^{tX} ) \big)|_{t=0}. $$
It is readily seen that 
\begin{equation} 
(\rho(X) f) (A) = df_A (AX) = \nabla f (A) \cdot AX = AX (f). 
\label{eq:rhoX_prop}
\end{equation}
In particular, $\rho$ is linear with respect to $X$. The map $\rho$: $\mathfrak{so}_n{\mathbb R} \to {\mathcal L} \big(C^\infty(\textrm{SO}_n{\mathbb R}) \big)$, $X \mapsto \rho(X)$ is actually a Lie algebra representation, i.e. on top of linearity, it satisfies 
$$ \rho([X,Y]) f = [\rho(X), \rho(Y)] f , \quad \forall X, Y \in \mathfrak{so}_n{\mathbb R}, \quad \forall f \in C^\infty(\textrm{SO}_n{\mathbb R}), $$
where ${\mathcal L} \big(C^\infty(\textrm{SO}_n{\mathbb R}) \big)$ is the space of linear maps from $C^\infty(\textrm{SO}_n{\mathbb R})$ into itself. Let us now introduce $F_{ij} =E_{ij} - E_{ji} = e_i \wedge e_j$, where the wedge product of two vectors $a$ and $b$ is their anti-symmetrized tensor product: $a \wedge b = a \otimes b - b \otimes a$. Then, $(F_{ij})_{1 \leq i < j \leq n}$ is an orthonormal basis of $\mathfrak{so}_n{\mathbb R}$ for the inner product~\eqref{eq:inner_product}. The Lie group Laplacian $\Delta_G$ is defined for a function $f \in C^\infty(\textrm{SO}_n{\mathbb R})$ as follows: 
\begin{equation} 
(\Delta_G f)(A)  = \sum_{1 \leq i<j \leq n} \frac{d^2}{dt^2} f(A e^{t F_{ij}}) \Big|_{t=0} = \sum_{1 \leq i<j \leq n} \big( \rho(F_{ij})^2 f \big) (A). 
\label{eq:laplie_def}
\end{equation}
It can be shown that the orthonormal basis $(F_{ij})_{1 \leq i < j \leq n}$ can be replaced by any other orthonormal basis of $\mathfrak{so}_n{\mathbb R}$ without changing the operator $\Delta_G$. 

Denote by ``$\textrm{Ad}$'' and ``$\textrm{ad}$'' the adjoint representations of $\textrm{SO}_n{\mathbb R}$ and $\mathfrak{so}_n{\mathbb R}$, respectively. $\textrm{Ad}$ maps $\textrm{SO}_n{\mathbb R}$ into the group $\textrm{Aut}(\mathfrak{so}_n{\mathbb R})$ of automorphisms of $\mathfrak{so}_n{\mathbb R}$ while $\textrm{ad}$ maps linearly $\mathfrak{so}_n{\mathbb R}$ into the space ${\mathcal L}(\mathfrak{so}_n{\mathbb R})$ of endomorphisms of $\mathfrak{so}_n{\mathbb R}$. They are defined for $A \in \textrm{SO}_n{\mathbb R}$ and $X, \, Y \in \mathfrak{so}_n{\mathbb R}$ by
$$ \textrm{Ad}(A) (Y) = A Y A^{-1},  \quad  \textrm{ad}(X) (Y) = [X,Y] = XY-YX. $$
$\textrm{Ad}$ is a Lie-group representation of $\textrm{SO}_n{\mathbb R}$ while $\textrm{ad}$ is a Lie-algebra representation of $\mathfrak{so}_n{\mathbb R}$, i.e. they satisfy
$$ \textrm{Ad}(A) \textrm{Ad}(B) = \textrm{Ad}(AB), \quad [\textrm{ad}(X), \textrm{ad}(Y)] = \textrm{ad}([X,Y]), $$
for all $A$, $B \in \textrm{SO}_n{\mathbb R}$ and all $X$, $Y \in \mathfrak{so}_n{\mathbb R}$. The latter formula is a re-expression of the Jacobi identity of the Lie bracket using its skew-commutativity. Now, for $A \in \textrm{SO}_n{\mathbb R}$, we have $A^{-1} = A^T$, and thus the inner product on $\mathfrak{so}_n{\mathbb R}$ satisfies
\begin{equation} 
\textrm{Ad}(A) (X) \cdot \textrm{Ad}(A) (Y) = X \cdot Y, \quad \textrm{ad}(Z) (X) \cdot Y + X \cdot \textrm{ad}(Z) (Y) = 0. 
\label{eq:invariant}
\end{equation}
An inner product satisfying \eqref{eq:invariant} is said to be invariant. Furthermore, denote by $L$ and $R$ the left and right regular representations of $\textrm{SO}_n{\mathbb R}$ on $C^\infty(\textrm{SO}_n{\mathbb R})$ defined for any $A, \, B \in \textrm{SO}_n{\mathbb R}$ and any $f \in C^\infty(\textrm{SO}_n{\mathbb R})$ by $(L(B) f)(A) = f(B^{-1} A)$ and $(R(B) f)(A) = f(AB)$. We have $L(B) (\Delta_G f) = \Delta_G (L(B) f)$, which means that the Lie-group Laplacian is left-invariant. Because the inner product on $\mathfrak{so}_n{\mathbb R}$ is invariant (see \eqref{eq:invariant}), we have  $R(B) (\Delta_G f) = \Delta_G (R(B) f)$, which means that the Lie-group Laplacian is also right-invariant.

It turns out that the two concepts of Laplacians coincide on $\textrm{SO}_n{\mathbb R}$ as stated in the following lemma, proved in Appendix \ref{secapp_lem:laplacians=_proof}. 

\begin{lemma}
We have $\Delta_M = \Delta_G$. 
\label{lem:laplacians=}
\end{lemma}

As a consequence, we can drop the subscript $M$ or $G$ of $\Delta$. 

Two elements $A$ and $B$ of $\textrm{SO}_n{\mathbb R}$ are said to be conjugate if and only if there exists $g \in \textrm{SO}_n{\mathbb R}$ such that $B = g A g^{-1}$. Conjugation is an equivalence relation and we write $A \sim B$ whenever $A$ and $B$ are conjugate. A class function is a function on $\textrm{SO}_n{\mathbb R}$ which is invariant by conjugation, i.e. which satifies:
$$ f(g A g^{-1}) = f(A), \quad \forall A, \, g \in \textrm{SO}_n{\mathbb R}. $$
A class function thus only depends on the conjugation classes. On $\textrm{SO}_n{\mathbb R}$ conjugation classes are well-known thanks to classical results on matrix reduction. Let $\theta \in {\mathbb R}$. We define the planar rotation matrix $R_\theta$ by
$$ R_\theta = \left( \begin{array}{rr} \cos \theta & - \sin \theta \\ \sin  \theta &  \cos \theta \end{array} \right). $$
For $(\theta_1, \ldots , \theta_p) \in {\mathbb R}^p$, $A_{\theta_1, \ldots ,\theta_p}$ denotes the following matrices defined by blocks:
\begin{itemize}
\item[-] in the case $n = 2p$, $p \geq 2$, 
\begin{equation} A_{\theta_1, \ldots , \theta_p} = \left( \begin{array}{cccc} 
\scalebox{1.2}{$R_{\theta_1}$} &  &  & \scalebox{2.}{$0$} \\
& \scalebox{1.2}{$R_{\theta_2}$} &  & \\
&  &  \ddots  &  \\
\scalebox{2.}{$0$} &  &  &  \scalebox{1.2}{$R_{\theta_p}$}
\end{array} \right) \in \mathrm{SO}_{2p}{\mathbb R}, 
\label{eq:R2p}
\end{equation}
\item[-] in the case $n = 2p+1$, $p \geq 1$, 
\begin{equation}
A_{\theta_1 , \ldots , \theta_p} = \left( \begin{array}{ccccc} 
\scalebox{1.2}{$R_{\theta_1}$} &  &  & \scalebox{2.}{$0$} & 0 \\
& \scalebox{1.2}{$R_{\theta_2}$} &  & & \vdots \\
&  &  \ddots  &  & \vdots \\
\scalebox{2.}{$0$} &  &  &  \scalebox{1.2}{$R_{\theta_p}$} & 0 \\
0  &  \ldots & \ldots & 0 & 1 
\end{array} \right) \in \mathrm{SO}_{2p+1}{\mathbb R}. 
\label{eq:R2p+1}
\end{equation}
\end{itemize}

We write $\Theta = (\theta_1,\theta_2, \ldots, \theta_p) \in {\mathbb R}^p$. We define the subset ${\mathbb T}$ of $\mathrm{SO}_{2p}{\mathbb R}$ or $\mathrm{SO}_{2p+1}{\mathbb R}$ by 
$$ {\mathbb T} = \{ A_\Theta  \, \, | \, \, \Theta \in {\mathbb R}^p \} .$$
${\mathbb T}$ is an abelian subgroup of $\mathrm{SO}_{2p}{\mathbb R}$ or $\mathrm{SO}_{2p+1}{\mathbb R}$. Consider the map $\psi$: $\Theta \to A_\Theta$. Since we have $A_\Theta A_{\Theta'} = A_{\Theta + \Theta'}$, $A_0 = \textrm{I}$, $\psi$ is a surjective group morphism ${\mathbb R}^p \to {\mathbb T}$ whose kernel is $(2 \pi {\mathbb Z})^p$. Thus, it induces a group isomorphism $\tilde \psi$, ${\mathcal T} =:  ({\mathbb R}/(2 \pi {\mathbb Z}))^p \to {\mathbb T}$. We note that ${\mathcal T} \cong ({\mathbb S}^1)^p$ is a $p$-dimensional torus. It can be shown that ${\mathbb T}$ is maximal among abelian subgoups of $\mathrm{SO}_n{\mathbb R}$. For that reason, it is called a maximal torus. 

Then, by matrix reduction theory, any element of $\mathrm{SO}_{2p}{\mathbb R}$ or $\mathrm{SO}_{2p+1}{\mathbb R}$ is conjugate to an element of ${\mathbb T}$.  Let $f$ be a $C^\infty$ class function on $\textrm{SO}_{2p} {\mathbb R}$ or $\textrm{SO}_{2p+1} {\mathbb R}$. Defining $\varphi_f(\Theta) = f(A_{\Theta})$, we have 
\begin{equation}
f(A) = \varphi_f(\Theta), \quad \forall A \, \textrm{ such that } \, A \sim A_{\Theta}, 
\label{eq:def_phif}
\end{equation}
and $\varphi_f \in C^\infty ({\mathcal T})$. Thus, the knowledge of a class function only requires the knowledge of a function of $\Theta \in {\mathcal T}$. This involves a dramatic reduction of complexity. 

However, elements of ${\mathbb T}$ which are conjugate to a given $A$ are not unique. The list of such elements differs in the even and odd dimensional cases: 
\begin{itemize}
\item[-] In the case of $\mathrm{SO}_{2p}{\mathbb R}$, $(\theta_1, \ldots, \theta_p)$ can undergo any permutations of $\{1, \ldots, p\}$. Permutation of $\theta_i$ and $\theta_j$ involves conjugation by the permutation matrix which exchanges $e_{2i-1}$ and $e_{2j-1}$ on the one hand and $e_{2i}$ and $e_{2j}$ on the other hand. Since the corresponding permutation is the product of two transpositions, it has signature $+1$ and so, the associated permutation matrix is an element of $\textrm{SO}_{2p} {\mathbb R}$. But one can also exchange $e_{2i-1}$ and $e_{2i}$ on the one hand, and $e_{2j-1}$ and $e_{2j}$ on the other hand. Again the corresponding permutation matrix is an element of $\textrm{SO}_{2p} {\mathbb R}$ and conjugating $A_\Theta$ by it amounts to making the change $(\ldots, \theta_i, \ldots, \theta_j, \ldots) \to (\ldots, -\theta_i, \ldots, -\theta_j, \ldots)$. Note that we couldn't just conjugate with the permutation matrix associated with the single transposition which exchanges $e_{2i-1}$ and $e_{2i}$ because the associated matrix is of determinant $-1$ and does not belong to $\textrm{SO}_{2p} {\mathbb R}$. Thus, only changes of an even number of signs in the $\theta_i$'s is allowed. We denote by ${\mathcal G}$ the group of bijections ${\mathcal T} \to {\mathcal T}$ generated by the transformations $(\ldots, \theta_i, \ldots, \theta_j, \ldots) \to (\ldots, \theta_j, \ldots, \theta_i, \ldots)$ and $(\ldots, \theta_i, \ldots, \theta_j, \ldots) \to (\ldots, -\theta_i, \ldots, -\theta_j, \ldots)$ where $(i,j)$ ranges over all pairs of integers in $\{1, \ldots, p\}$ such that $i \not = j$. 

\item[-] In the case of $\mathrm{SO}_{2p+1}{\mathbb R}$, $(\theta_1, \ldots, \theta_p)$ can again  undergo any permutations of $\{1, \ldots, p\}$ for the same reason as in the even-dimensional case. But one can also conjugate by a matrix which exchanges $e_{2i-1}$ and $e_{2i}$ and bears a $-1$ in the last diagonal element (which was not possible in the odd-dimensional case). This matrix has determinant unity and is thus an element of $\textrm{SO}_{2p+1} {\mathbb R}$. Conjugating $A_\Theta$ by it amounts to making the change $(\ldots, \theta_i, \ldots) \to (\ldots, -\theta_i, \ldots)$. Thus, changes of an arbitrary number of signs of the $\theta_i$'s are allowed in this case. In this case, ${\mathcal G}$ is generated by the transformations $(\ldots, \theta_i, \ldots, \theta_j, \ldots) \to (\ldots, \theta_j, \ldots, \theta_i, \ldots)$ where $(i,j)$ ranges over all pairs of integers in $\{1, \ldots, p\}$ such that $i \not = j$ and by $(\ldots, \theta_i, \ldots) \to (\ldots, -\theta_i, \ldots)$, where $i$ ranges over all integers in $\{1, \ldots, p\}$. 
\end{itemize}

In both the even and odd dimensional cases, all conjugations mapping two elements of ${\mathbb T}$ are generated by the conjugations described above. To see that we introduce the set of all the elements of $\textrm{SO}_n {\mathbb R}$ which conjugate an element of ${\mathbb T}$ into an element of ${\mathbb T}$. This set is a group called the normalizer of ${\mathbb T}$ and is defined by 
$$ N({\mathbb T}) = \big\{  g \in \textrm{SO}_n {\mathbb R} \, \, | \, \, g {\mathbb T} g^{-1} = {\mathbb T} \big\}. $$
Let us introduce the group homomorphism $\Phi$: $N({\mathbb T}) \to \textrm{Aut}({\mathbb T})$, $g \mapsto \Phi_g$ (where $\textrm{Aut}({\mathbb T})$ denotes the group of group automorphisms of ${\mathbb T}$), such that $\Phi_g(h) = ghg^{-1}$, $\forall h \in {\mathbb T}$. This homomorphism has kernel equal to ${\mathbb T}$ (that ${\mathbb T} \subset \textrm{ker} \, \Phi$ is obvious. The converse is due to the fact that ${\mathbb T}$ is a maximal abelian subgroup of $\textrm{SO}_n {\mathbb R}$). Thus, $\Phi$ induces a group isomorphism $\tilde \Phi$: $N({\mathbb T})/{\mathbb T} \to \textrm{im} (\Phi) \subset \textrm{Aut}({\mathbb T})$. Now, define the isomorphism $\Psi$: $\textrm{Aut}({\mathbb T}) \to \textrm{Aut}({\mathcal T})$, $\varphi \mapsto \tilde \psi^{-1} \circ \varphi \circ \tilde \psi$, $\forall \varphi \in \textrm{Aut}({\mathbb T})$, where $\tilde \psi$: ${\mathcal T} \to {\mathbb T}$ was defined above. Then, $\Psi \circ \tilde \Phi$ is an isomorphism between $N({\mathbb T})/{\mathbb T}$ and the group of transformations of ${\mathcal T}$ that are induced by conjugations. It turns out that the group $N({\mathbb T})/{\mathbb T}$ is well-known to be isomorphic to the Weyl group $\mathfrak{W}$ \cite[Claim 26.15]{fulton2013representation}. $\mathfrak{W}$ is a finite subsgroup of the orthogonal group of the space ${\mathbb R}^p$. If ${\mathbb R}^p$ is endowed with the canonical basis $(e_i)_{i=1}^p$, then $\mathfrak{W}$ is generated by the orthogonal symmetries in the hyperplanes $\{e_i \pm e_j\}^\bot$, with $i$, $j \in \{1, \ldots, p\}$, $i \not = j$ in the case $n=2p$ and by the orthogonal symmetries in the hyperplanes $\{e_i \pm e_j\}^\bot$ and $\{e_i\}^\bot$ in the case $n=2p+1$ \cite[p. 271]{fulton2013representation}. But it is readily checked that the group ${\mathcal G}$ of transformations of ${\mathcal T}$ defined above is also generated by the same transformations. Hence, ${\mathcal G} = \mathfrak{W}$ showing the claim. 

It results that if $f \in C^\infty$ is a class function, the function $\varphi_f \in C^\infty({\mathcal T})$ associated with it through \eqref{eq:def_phif} must be invariant by the action of $\mathfrak{W}$. Hence, $\varphi_f$ is a function defined on the quotient ${\mathcal T}' =: {\mathcal T}/\mathfrak{W}$. For simplicity, we will drop the prime and from now on, ${\mathcal T}$ will actually mean ${\mathcal T}'$. 

Now, because $\Delta$ is left and right invariant, it is invariant by conjugation (indeed, conjugation by $g$ is the composition $L(g^{-1})R(g^{-1})$). It results that if $f$ is a class function, $\Delta f$ is a class function. For any $\varphi \in C^{\infty}({\mathcal T})$, there exists a unique function $f \in C^\infty(\textrm{SO}_n{\mathbb R})$ such that $\varphi = \varphi_f$. Indeed, $f$ is given by $f(A) = \varphi(\Theta)$, for any $\Theta$ such that $A \sim A_\Theta$. Then, we can define a linear operator $L$: $C^{\infty}({\mathcal T}) \to C^{\infty}({\mathcal T})$ by 
$$ L \varphi_f = \varphi_{\Delta f}, \quad \forall f \in C^{\infty}({\mathbb T}). $$
The operator $L$ is called the radial Laplacian. The goal of this paper is to find an explicit expression of $L$ and to state some of its properties.

\setcounter{equation}{0}
\section{Radial Laplacian}
\label{sec_rad_lap}

\subsection{Expression}
\label{subsec_statement}

The expression of the radial Laplacian is given by the following: 

\begin{theorem}
The radial Laplacian on the rotation groups is given by the following expressions: 

\noindent
(i) Case of $\mathrm{SO}_{2p}{\mathbb R}$ for $p \geq 2$: 
\begin{equation} L =  \sum_{j=1}^p \frac{\partial^2}{\partial \theta_j^2} + \sum_{1 \leq j < k \leq p} \frac{2}{\cos \theta_k - \cos \theta_j} \, \Big( \big( \sin \theta_j \, \frac{\partial}{\partial \theta_j} - \sin \theta_k \, \frac{\partial}{\partial \theta_k}) \Big) ; 
\label{eq:radlap_2p}
\end{equation}

\noindent
(ii) Case of $\mathrm{SO}_{2p+1}{\mathbb R}$ for $p \geq 1$: 
\begin{eqnarray}
&&\hspace{-1cm}
 L =  \sum_{j=1}^p \Big(\frac{\partial^2}{\partial \theta_j^2} + \frac{\sin \theta_j}{1 - \cos \theta_j} \, \frac{\partial}{\partial \theta_j} \Big) \nonumber \\
&&\hspace{1cm}
+ \sum_{1 \leq j < k \leq p} \frac{2}{\cos \theta_k - \cos \theta_j} \, \Big( \big( \sin \theta_j \, \frac{\partial}{\partial \theta_j} - \sin \theta_k \, \frac{\partial}{\partial \theta_k}) \Big) . 
\label{eq:radlap_2p+1}
\end{eqnarray}
We note that in the case $p=1$, the second sum simply disappears in \eqref{eq:radlap_2p+1}. We leave aside the case $\mathrm{SO}_{2}{\mathbb R}$ which is trivial, the group being abelian. 
\label{thm:statement}
\end{theorem}

To prove this theorem, we will first review the cases of $\textrm{SO}_3 {\mathbb R}$ and $\textrm{SO}_4 {\mathbb R}$ and then, generalize the analysis to the general cases $\mathrm{SO}_{2p}{\mathbb R}$ and $\mathrm{SO}_{2p+1}{\mathbb R}$. We will adopt the method presented in \cite{faraud2008Analysis} (cf Proof of Proposition 8.3.3 and Section 10.2). This method relies on the following lemma: 

\begin{lemma}
If $f$ is a class function on $\textrm{SO}_n {\mathbb R}$, for all $X \in \mathfrak{so}_n {\mathbb R}$ and all $A \in \textrm{SO}_n {\mathbb R}$, we have
\begin{equation}
\rho \big( \textrm{Ad}(A^{-1}) X -X \big)^2 f(A) - \rho \big( [\textrm{Ad}(A^{-1}) X, X ] \big) f(A) = 0. 
\label{eq:fonda}
\end{equation}
\label{lem:fonda}
\end{lemma}

\noindent
\textbf{Proof.} If $f$ is a class function, we have $f(A) = f(e^{tX} A e^{-tX})$, for all $A \in \textrm{SO}_n {\mathbb R}$ and all $X \in \mathfrak{so}_n{\mathbb R}$. This formula can be recast into $f(A) = f(A e^{t \textrm{Ad}(A^{-1}) X} e^{-tX})$. Taking the second derivative of this formula with respect to $t$ at $t=0$ and using that 
\begin{equation}
\frac{d^2}{dt^2} f(A e^{tX} e^{tY}) \big|_{t=0} = \rho(X+Y)^2 f(A) + \rho \big( [X,Y] \big) f (A), \quad \forall X, \, Y \in \mathfrak{so}_n{\mathbb R}, 
\label{eq:secder}
\end{equation}
we get \eqref{eq:fonda}. Eq. \eqref{eq:secder} is proved in \cite[Section 8.2 (d)]{faraud2008Analysis}. \endproof

We will use the following formula, whose proof is straightforwad: 
\begin{equation}
[F_{ij},F_{k \ell}] = \big( \delta_{jk} F_{i \ell} + \delta_{i \ell} F_{jk} - \delta_{ik} F_{j \ell} - \delta_{j \ell} F_{ik} \big). 
\label{eq:commutF}
\end{equation}
We will also use the notation $\{A,B\}$ for the anti-commutator of two operators $A$, $B$, i.e. 
$$ \{A, B \} = AB + BA. $$

\subsection{Proof of Theorem \ref{thm:statement}: case of $\textrm{SO}_3 {\mathbb R}$}
\label{subsec_SO3}

This is the case $\textrm{SO}_{2p+1} {\mathbb R}$ with $p=1$. Thus, the matrix \eqref{eq:R2p+1} has a single block $R_\theta$. So, the vector $\Theta$ has a single coordinate, $\theta$ and we use the notation $A_\theta$ for $A_\Theta$. Let $f$ be a $C^\infty$ class function on $\textrm{SO}_3 {\mathbb R}$ and $\varphi_f = \varphi_f(\theta)$ the corresponding function of $\theta$.  

Define
\begin{eqnarray*} G^+ &=& \frac{1}{\sqrt{2}} (F_{13} + F_{23}) = \frac{1}{\sqrt{2}} (e_1 \otimes e_3 + e_2 \otimes e_3 - e_3 \otimes e_1 - e_3 \otimes e_2), \\
G^- &=& \frac{1}{\sqrt{2}} (F_{13} - F_{23}) = \frac{1}{\sqrt{2}} (e_1 \otimes e_3 - e_2 \otimes e_3 - e_3 \otimes e_1 + e_3 \otimes e_2). 
\end{eqnarray*}
Then, $(F_{12}, G^+, G^-)$ is still an orthonormal basis of $\mathfrak{so}_3 {\mathbb R}$. We choose this basis to express the Laplacian of a $C^\infty$ function $f$ on $\textrm{SO}_3 {\mathbb R}$: 
\begin{equation} \Delta f(A) = \big( \rho(F_{12})^2 f \big)(A) + \big( \rho(G^+)^2 f \big)(A)  + \big( \rho(G^-)^2 f \big)(A) . 
\label{eq:Delta_SO3}
\end{equation}

We now use \eqref{eq:fonda} with $A=A_\theta$ and $X = G^+$ and $X = G^-$ successively. We have
$$\textrm{Ad}(A_{-\theta}) G^+ = \cos \theta \, G^+ + \sin \theta \, G^- , \quad
\textrm{Ad}(A_{-\theta}) G^- = - \sin \theta \, G^+ + \cos \theta \, G^- .  
$$
Then, using the linearity of the map $X \mapsto (\rho(X) f) (A)$ for any $A \in \textrm{SO}_2{\mathbb R}$, we get:
\begin{eqnarray*} 
\rho \big( \textrm{Ad}(A_{-\theta}) G^+ - G^+ \big)^2 f (A_\theta) &=& \Big( \big( (\cos \theta - 1)^2 \, \rho(G^+)^2 + \sin^2 \theta \, \rho(G^-)^2 \\
& & \hspace{1cm} + (\cos \theta - 1) \, \sin \theta \, \big\{ \rho(G^+), \rho(G^-) \big\} \big) f \Big) (A_\theta),  \\
\rho \big( \textrm{Ad}(A_{-\theta}) G^- - G^- \big)^2 f (A_\theta) &=& \Big( \big( (\cos \theta - 1)^2 \, \rho(G^-)^2 + \sin^2 \theta \, \rho(G^+)^2 \\
& & \hspace{1cm} - (\cos \theta - 1) \, \sin \theta \, \big\{ \rho(G^+), \rho(G^-) \big\}  \big) f \Big) (A_\theta). 
\end{eqnarray*}
Furthermore, using \eqref{eq:commutF}, we get
\begin{eqnarray*} [\textrm{Ad}(A_{-\theta}) G^+, G^+] &=& [\textrm{Ad}(A_{-\theta}) G^-, G^-] = \sin \theta \, [G^-, G^+] \\
&=& \sin \theta \, [F_{13}, F_{23}] = - \sin \theta \, F_{12}. 
\end{eqnarray*}
Thus, adding the equations \eqref{eq:fonda} corresponding to $G^+ $ and $G^-$ together leads to 
\begin{equation} 
2 (1 - \cos \theta) \Big( \big( \rho(G^+)^2 + \rho(G^-)^2 \big) f \Big) (A_\theta) = - 2 \, \sin \theta \,  \big( \rho(F_{12}) f \big)(A_\theta). 
\label{eq:rhoGpmsq}
\end{equation}
Thus, with \eqref{eq:Delta_SO3}, we get
\begin{equation} 
(\Delta f)(A_\theta) = \Big( \big( \rho(F_{12})^2 - \frac{\sin \theta}{1 - \cos \theta} \, \rho(F_{12}) \big) f \Big) (A_\theta). 
\label{eq:radlap_SO3_2}
\end{equation}
Now, we remark that $e^{t F_{12}} = A_{-t}$. Thus, by the definition of $\rho$: 
\begin{eqnarray} 
\big( \rho(F_{12}) f \big) (A_\theta) &=& \frac{d}{dt} f(A_\theta e^{t F_{12}}) \big|_{t=0} = \frac{d}{dt} f(A_{\theta - t}) \big|_{t=0} \\
&=& \frac{d}{dt} \varphi_f(\theta - t) \big|_{t=0} = - \frac{d}{d\theta} \varphi_f (\theta). 
\label{eq:F12=dtheta}
\end{eqnarray}
Inserting this expresion into \eqref{eq:radlap_SO3_2}, we get the expression of the radial Laplacian on $\textrm{SO}_3{\mathbb R}$: 
$$ L =  \big( \frac{d^2}{d \theta^2} + \frac{\sin \theta}{1 - \cos \theta} \, \frac{d}{d\theta} \big), $$
which corresponds to \eqref{eq:radlap_2p+1} in the case $p=1$. \endproof

\subsection{Proof of Theorem \ref{thm:statement}: case of $\textrm{SO}_4 {\mathbb R}$}
\label{subsec_SO4}

This is the case $\textrm{SO}_{2p} {\mathbb R}$ with $p=2$. Thus, the matrix \eqref{eq:R2p} has a two blocks $R_{\theta_1}$ and $R_{\theta_2}$. So, $\Theta=(\theta_1,\theta_2)$. Let $f$ be a $C^\infty$ class function on $\textrm{SO}_4 {\mathbb R}$ and $\varphi_f = \varphi_f(\theta_1,\theta_2)$ the corresponding function of $(\theta_1,\theta_2)$.  

We now define
\begin{eqnarray*} H^+ &=& \frac{1}{\sqrt{2}} (F_{13} + F_{24}), \quad H^- = \frac{1}{\sqrt{2}} (F_{13} - F_{24}) , \\
K^+ &=& \frac{1}{\sqrt{2}} (F_{14} + F_{23}), \quad K^- = \frac{1}{\sqrt{2}} (F_{14} - F_{23}) 
\end{eqnarray*}
Then, $(F_{12}, F_{34}, H^+, H^-, K^+, K^-)$ is still an orthonormal basis of $\mathfrak{so}_4 {\mathbb R}$. We choose this basis to express the Laplacian of a $C^\infty$ function $f$ on $\textrm{SO}_4 {\mathbb R}$: 
\begin{equation} \Delta f(A) = \Big( \big( \rho(F_{12})^2 + \rho(F_{34})^2 + \rho(H^+)^2 + \rho(H^-)^2 + \rho(K^+)^2 + \rho(K^-)^2 \big) f \Big)(A). 
\label{eq:Delta_SO4}
\end{equation}

We now use \eqref{eq:fonda} with $A=A_{\theta_1,\theta_2}$ and $X = H^+$, $H^-$, $K^+$, $K^-$ successively. We have
\begin{eqnarray*}
\textrm{Ad}(A_{-\theta_1,-\theta_2}) H^+ &=& (c_1 c_2 + s_1 s_2) H^+ - (c_1 s_2- s_1 c_2) K^- , \\
\textrm{Ad}(A_{-\theta_1,-\theta_2}) H^- &=& (c_1 c_2 - s_1 s_2) H^- - (c_1 s_2 + s_1 c_2) K^+ , \\
\textrm{Ad}(A_{-\theta_1,-\theta_2}) K^+ &=& (c_1 s_2 + s_1 c_2) H^- +  (c_1 c_2 - s_1 s_2) K^+ , \\
\textrm{Ad}(A_{-\theta_1,-\theta_2}) K^- &=& (c_1 s_2- s_1 c_2) H^+ + (c_1 c_2 + s_1 s_2) K^-, 
\end{eqnarray*} 
with $c_i = \cos \theta_i$, $s_i = \sin \theta_i$ for $i=1, \, 2$. Then, we get:
\begin{eqnarray*} 
&&\hspace{-1cm} 
\rho \big( \textrm{Ad}(A_{-\theta_1,-\theta_2}) H^+ - H^+ \big)^2 f (A_{\theta_1,\theta_2}) = \Big( \big( (c_1 c_2 + s_1 s_2 - 1)^2 \, \rho(H^+)^2 \\
& & \hspace{2cm} 
+ (c_1 s_2- s_1 c_2)^2 \, \rho(K^-)^2 \\
& & \hspace{2cm} 
- (c_1 c_2 + s_1 s_2 - 1) \, (c_1 s_2- s_1 c_2) \, \big\{ \rho(H^+), \rho(K^-) \big\} \big) f \Big) (A_{\theta_1,\theta_2}),  \\
&&\hspace{-1cm}
\rho \big( \textrm{Ad}(A_{-\theta_1,-\theta_2}) H^- - H^- \big)^2 f (A_{\theta_1,\theta_2}) = \Big( \big( (c_1 c_2 - s_1 s_2 - 1)^2 \, \rho(H^-)^2 \\
& & \hspace{2cm} 
+ (c_1 s_2 + s_1 c_2)^2 \, \rho(K^+)^2 \\
& & \hspace{2cm} 
- (c_1  c_2 - s_1 s_2 - 1) \, (c_1 s_2 + s_1 c_2) \, \big\{ \rho(H^-), \rho(K^+) \big\} \big) f \Big) (A_{\theta_1,\theta_2}),  
\end{eqnarray*}

\begin{eqnarray*}
&&\hspace{-1cm}
\rho \big( \textrm{Ad}(A_{-\theta_1,-\theta_2}) K^+ - K^+ \big)^2 f (A_{\theta_1,\theta_2}) = \Big( \big( (c_1 s_2 + s_1 c_2)^2 \, \rho(H^-)^2 \\
& & \hspace{2cm} 
+ (c_1 c_2 - s_1 s_2 - 1)^2 \, \rho(K^+)^2 \\
& & \hspace{2cm} 
+ (c_1 s_2 + s_1 c_2) \, (c_1 c_2 - s_1 s_2 - 1) \, \big\{ \rho(H^-), \rho(K^+) \big\} \big) f \Big) (A_{\theta_1,\theta_2}),  \\
&&\hspace{-1cm}
\rho \big( \textrm{Ad}(A_{-\theta_1,-\theta_2}) K^- - K^- \big)^2 f (A_{\theta_1,\theta_2}) = \Big( \big( (c_1 s_2- s_1 c_2)^2 \, \rho(H^+)^2  \\
& & \hspace{2cm}
+ (c_1 c_2 + s_1 s_2 - 1)^2 \, \rho(K^-)^2 \\
& & \hspace{2cm}
+ (c_1 s_2- s_1 c_2) \, (c_1 c_2 + s_1 s_2 - 1) \, \big\{ \rho(H^+), \rho(K^-) \big\} \big) f \Big) (A_{\theta_1,\theta_2}). 
\end{eqnarray*}
Furthermore, using \eqref{eq:commutF}, we get
\begin{eqnarray*} 
&& \hspace{-1cm} 
[\textrm{Ad}(A_{-\theta_1,-\theta_2}) H^+, H^+] = [\textrm{Ad}(A_{-\theta_1,-\theta_2}) K^-, K^-] = - (c_1 s_2- s_1 c_2) [K^-,H^+] \\
&& \hspace{7cm} 
=  (c_1 s_2- s_1 c_2) (F_{12} - F_{34}), \\
&& \hspace{-1cm} 
\mbox{}[\textrm{Ad}(A_{-\theta_1,-\theta_2}) H^-, H^-] = [\textrm{Ad}(A_{-\theta_1,-\theta_2}) K^+, K^+] = - (c_1 s_2 + s_1 c_2) [K^+,H^-] \\
&& \hspace{7cm} 
= - (c_1 s_2 + s_1 c_2) (F_{12} + F_{34}), 
\end{eqnarray*}
Thus, adding the equations \eqref{eq:fonda} corresponding to $H^+ $ and $K^-$ (resp. $H^- $ and $K^+$) together leads to 
\begin{eqnarray*} 
\Big( \big( \rho(H^+)^2 + \rho(K^-)^2 \big) f \Big) (A_{\theta_1,\theta_2}) &=&  \frac{\sin (\theta_2 - \theta_1)}{1 - \cos (\theta_2 - \theta_1)} \,  \Big( \big( \rho(F_{12}) - \rho(F_{34}) \big) f \Big)(A_{\theta_1,\theta_2}) , \\
\Big( \big( \rho(H^-)^2 + \rho(K^+)^2 \big) f \Big) (A_{\theta_1,\theta_2}) &=& - \, \frac{\sin (\theta_1 + \theta_2)}{1 - \cos (\theta_1 + \theta_2)} \,  \Big( \big( \rho(F_{12}) + \rho(F_{34}) \big) f \Big)(A_{\theta_1,\theta_2}) , 
\end{eqnarray*}
Thus, with \eqref{eq:Delta_SO4}, we get
\begin{eqnarray} 
&& \hspace{-1cm} 
(\Delta f)(A_{\theta_1,\theta_2}) = \Big( \Big( \rho(F_{12})^2 + \rho(F_{34})^2 +  \frac{\sin (\theta_2 - \theta_1)}{1 - \cos (\theta_2 - \theta_1)} \,  \big( \rho(F_{12}) - \rho(F_{34}) \big) \nonumber \\
&& \hspace{3cm} 
- \frac{\sin (\theta_1 + \theta_2)}{1 - \cos (\theta_1 + \theta_2)} \,  \big( \rho(F_{12}) + \rho(F_{34}) \big) \Big) f \Big) (A_{\theta_1,\theta_2}). 
\label{eq:radlap_SO4_2}
\end{eqnarray}
Now, similarly to \eqref{eq:F12=dtheta}, we have 
$$ \big( \rho(F_{12}) f \big)(A_{\theta_1,\theta_2}) = -  \frac{\partial \varphi_f}{\partial \theta_1} (\theta_1,\theta_2), \quad \big( \rho(F_{34}) f \big)(A_{\theta_1,\theta_2}) = -  \frac{\partial \varphi_f}{\partial \theta_2} (\theta_1,\theta_2). $$
Inserting it into \eqref{eq:radlap_SO4_2} leads to 
\begin{eqnarray*}
L &=& \frac{\partial^2}{\partial \theta_1^2} + \frac{\partial^2}{\partial \theta_2^2} + \frac{\sin (\theta_1 - \theta_2)}{1 - \cos (\theta_1 - \theta_2)} \,  \big( \frac{\partial}{\partial \theta_1} - \frac{\partial}{\partial \theta_2} \big) \\
&& \hspace{5cm} 
+  \frac{\sin (\theta_1 + \theta_2)}{1 - \cos (\theta_1 + \theta_2)} \,  \big( \frac{\partial}{\partial \theta_1} + \frac{\partial}{\partial \theta_2} \big).   
\end{eqnarray*}
This can be recast into 
$$L = \frac{\partial^2}{\partial \theta_1^2} + \frac{\partial^2}{\partial \theta_2^2} + \frac{2}{\cos \theta_2 - \cos \theta_1} \, \big( \sin \theta_1 \, \frac{\partial}{\partial \theta_1} - \sin \theta_2 \, \frac{\partial}{\partial \theta_2} \big),   $$
which corresponds to \eqref{eq:radlap_2p} in the case $p=1$.  \endproof

\subsection{Proof of Theorem \ref{thm:statement}: case of $\textrm{SO}_{2p} {\mathbb R}$, $p \geq 1$}
\label{subsec_SO2p}

For $1 \leq j < k \leq p$, we define
\begin{eqnarray} H_{jk}^+ &=& \frac{1}{\sqrt{2}} (F_{2j-1 \, 2k-1} + F_{2j \, 2k}), \quad H_{jk}^- = \frac{1}{\sqrt{2}} (F_{2j-1 \, 2k-1} - F_{2j \, 2k}) , \label{eq:Hij}\\
K_{jk}^+ &=& \frac{1}{\sqrt{2}} (F_{2j-1 \, 2k} + F_{2j \, 2k-1}), \quad K_{jk}^- = \frac{1}{\sqrt{2}} (F_{2j-1 \, 2k} + F_{2j \, 2k-1}) . \label{eq:Kij}
\end{eqnarray}
Then, $\big( (F_{2j-1 \, 2j})_{j=1, \ldots, p}, (H_{jk}^+, H_{jk}^-, K_{jk}^+, K_{jk}^-)_{1 \leq j < k \leq p} \big)$ is still an orthonormal basis of $\mathfrak{so}_{2p} {\mathbb R}$. We choose this basis to express the Laplacian of a $C^\infty$ function $f$ on $\textrm{SO}_{2p} {\mathbb R}$: 
\begin{equation} \Delta f(A) = \Big( \Big( \sum_{j=1}^p \rho(F_{j \, j+1})^2 + \sum_{1 \leq j < k \leq p} \big( \rho(H_{jk}^+)^2 + \rho(H_{jk}^-)^2 + \rho(K_{jk}^+)^2 + \rho(K_{jk}^-)^2 \big) \Big) f \Big)(A). 
\label{eq:Delta_SO2p}
\end{equation}

We note that the computations of $\textrm{Ad}(A_{-\Theta}) H_{jk}^\pm$ and $\textrm{Ad}(A_{- \Theta}) K_{jk}^\pm$ only involve the four $2 \times 2$ matrix subblocks corresponding to line and column indices belonging to $\{ 2j-1, 2j \} \times \{ 2j-1, 2j \}$, $\{ 2j-1, 2j \} \times \{ 2k-1, 2k \}$, $\{ 2k-1, 2k \} \times \{ 2j-1, 2j \}$ and $\{ 2k-1, 2k \} \times \{ 2k-1, 2k \}$. Thus, the computation reduces to computations on $4 \times 4$ matrices that are identical to what was done in the case of $\textrm{SO}_{4}{\mathbb R}$. Thus, we can directly write that for any $1 \leq j < k \leq p$, we have  
\begin{eqnarray} 
&&\hspace{-1cm}
\Big( \big( \rho(H_{jk}^+)^2 + \rho(K_{jk}^-)^2 + \rho(H_{jk}^-)^2 + \rho(K_{jk}^+)^2\big) f \Big) (A_\Theta) \\
&&\hspace{-0.5cm}
= - \frac{2}{\cos \theta_k - \cos \theta_j} \, \Big( \big( \sin \theta_j \, \rho(F_{2j-1 \, 2j}) - \sin \theta_k \, \rho(F_{2k-1 \, 2k}) \big) f \Big) (A_\Theta).
\label{eq:rhoHKijsq}
\end{eqnarray}
Now, similarly to \eqref{eq:F12=dtheta}, for $j \in \{1, \ldots, p \}$, we have 
$$ \big( \rho(F_{2j-1 \, 2j}) f \big) (A_\Theta) = -  \frac{\partial \varphi_f}{\partial \theta_j} (\Theta). $$
Inserting these into \eqref{eq:Delta_SO2p}, we get \eqref{eq:radlap_2p}.  \endproof

\subsection{Proof of Theorem \ref{thm:statement}: case of $\textrm{SO}_{2p+1} {\mathbb R}$, $p \geq 1$}
\label{subsec_SO2p+1}

For $1 \leq j < k \leq p$, we define $H_{jk}^\pm$ and $K_{jk}^\pm$ in a similar way as in \eqref{eq:Hij}, \eqref{eq:Kij}. For $j \in \{1, \ldots, p\}$, we define 
$$G_j^+ = \frac{1}{\sqrt{2}} (F_{2j-1 \, 2p+1} + F_{2j \, 2p+1}), \quad
G_j^- = \frac{1}{\sqrt{2}} (F_{2j-1 \, 2p+1} - F_{2j \, 2p+1}). 
$$
Then, $\big( (F_{2j-1 \, 2j}, G_j^+, G_j^- )_{j=1, \ldots, p}, (H_{jk}^+, H_{jk}^-, K_{jk}^+, K_{jk}^-)_{1 \leq j < k \leq p} \big)$ is still an orthonormal basis of $\mathfrak{so}_{2p+1} {\mathbb R}$. We choose this basis to express the Laplacian of a $C^\infty$ function $f$ on $\textrm{SO}_{2p+1} {\mathbb R}$: 
\begin{eqnarray} \Delta f(A) &=& \Big( \Big( \sum_{j=1}^p \big( \rho(F_{j \, j+1})^2 + \rho(G_j^+)^2 + \rho(G_j^-)^2 \big)  \nonumber \\
&&\hspace{0cm}
+ \sum_{1 \leq j < k \leq p} \big( \rho(H_{jk}^+)^2 + \rho(H_{jk}^-)^2 + \rho(K_{jk}^+)^2 + \rho(K_{jk}^-)^2 \big) \Big) f \Big)(A). 
\label{eq:Delta_SO2p+1}
\end{eqnarray}
The computations involving $H_{jk}^\pm$ and $K_{jk}^\pm$ are identical to those of the preceding section. Those involving $G_j^\pm$ are identical to those done in Section \ref{subsec_SO3} for $\mathrm{SO}_3 {\mathbb R}$. Thus, we can collect \eqref{eq:rhoHKijsq} for $1 \leq j < k \leq p$ and \eqref{eq:rhoGpmsq} with $G^\pm$ replaced by $G_j^\pm$, $\theta$ by $\theta_j$ and $F_{12}$ by $F_{2j-1 \, 2j}$, for all $j \in \{1, \ldots, p \}$, and get \eqref{eq:radlap_2p+1}.  \endproof

\setcounter{equation}{0}
\section{Some properties of the radial Laplacian}
\label{sec_prop}

\subsection{Alternate expressions}
\label{subsec_alternate}

As a preliminary to this section, we state the Weyl integration formula \cite[Theorems IX.9.4 \& IX.9.5]{Simon}. We introduce the following functions: 
\begin{eqnarray}
\hspace{-0.5cm} \Pi_{2p} (\Theta) &=& \prod_{1 \leq j<k \leq p} \big( \cos \theta_j - \cos \theta_k \big), \quad \textrm{ for } p \geq 2, 
\label{eq:Pi2p} \\
\hspace{-0.5cm} \Pi_{2p+1}  (\Theta) &=&  \prod_{1 \leq j<k \leq p} \big( \cos \theta_j - \cos \theta_k \big) \, \prod_{j=1}^p \sin \frac{\theta_j}{2}, \quad \textrm{ for } p \geq 1.  
\label{eq:Pi2p+1}
\end{eqnarray}
We define 
\begin{equation}
\hspace{-0.5cm} u_{2p} = \frac{2^{(p-1)^2}}{p!} \Pi_{2p}^2, \, \,  \textrm{ for } p \geq 2; \qquad u_{2p+1} = \frac{2^{p(p-1)}}{p!} \Pi_{2p+1}^2, \, \,  \textrm{ for } p \geq 1. 
\label{eq:un}
\end{equation}
Then, we have

\begin{proposition}[Weyl integration formula]~

\noindent
Let $n \in {\mathbb N}$, $n\geq 3$. Let $p \in {\mathbb N}$ such that $n=2p$ or $n=2p+1$. For any integrable class function $f$ on $\mathrm{SO}_n{\mathbb R}$, we have 
\begin{equation}
\int_{\mathrm{SO}_n{\mathbb R}} f(A) \, dA = \frac{1}{(2 \pi)^p} \int_{[0,2 \pi]^p} f(A_{\Theta}) \, u_n (\Theta) \, d \Theta,
\label{eq:WIF}
\end{equation}
with $d \Theta = d \theta_1 \ldots d \theta_p$. 
\label{prop:WIF}
\end{proposition}

For a smooth function $\varphi$: ${\mathcal T} \to {\mathbb R}^p$ , $\nabla_\Theta \varphi$ denotes the vector field ${\mathcal T} \to {\mathbb R}^p$ of coordinates $(\frac{\partial \varphi}{\partial \theta_i})_{i=1}^p$. Similarly, for a smooth vector field ${\mathcal X}$: ${\mathcal T} \to {\mathbb R}^p$ of coordinates $({\mathcal X}_i)_{i=1}^p$, the divergence $\nabla_\Theta \cdot {\mathcal X}$ is the scalar $\sum_{i=1}^p \frac{\partial {\mathcal X}_i}{\partial \theta_i}$. The Laplacian $\Delta_\Theta$ is naturally defined by $\Delta_\Theta = \nabla_\Theta \cdot \nabla_\Theta \varphi$. Finally, we introduce the constants 
\begin{equation}
\gamma_{2p} = \sum_{j=1}^p (p-j)^2, \quad \gamma_{2p+1} = \sum_{j=1}^p (p-j+\frac{1}{2})^2. 
\label{eq:gamma}
\end{equation}
We note that $\gamma_{2p}$ (for $p \geq 2$) and $\gamma_{2p+1}$ (for $p \geq 1$) are strictly positive constants. We have the following 

\begin{lemma} (i) Case $n=2p$, $p \geq 2$: we have 
\begin{equation}
\nabla_\Theta u_{2p} (\Theta) = \Big( \big( \sum_{k \not = j} \frac{2}{\cos \theta_k - \cos \theta_j} \big) \sin \theta_j \,  u_{2p} (\Theta) \Big)_{j=1}^p.  
\label{eq:na_Th_u2p} 
\end{equation}

\smallskip
\noindent
(ii) Case $n=2p+1$, $p \geq 1$: we have 
\begin{equation}
\nabla_\Theta u_{2p+1} (\Theta) = \Big( \big( \sum_{k \not = j} \frac{2}{\cos \theta_k - \cos \theta_j} + \frac{1}{1 - \cos \theta_j} \big) \sin \theta_j  \,  u_{2p+1} (\Theta) \Big)_{j=1}^p.  
\label{eq:na_Th_u2p+1} 
\end{equation}

\smallskip
\noindent
(iii) In both the cases $n=2p$ ($p \geq 2$) and $n=2p+1$ ($p \geq 1$), we have 
\begin{equation}
\Delta_\Theta \Pi_n = - \gamma_n \Pi_n. 
\label{eq:Delt_Pin}
\end{equation}
\label{lem_auxi_un_Pin}
\end{lemma}

\noindent
\textbf{Proof.} (i) We compute
\begin{eqnarray*}
&&\hspace{-0.5cm}
\frac{\partial}{\partial \theta_j} \log u_{2p} = 2 \frac{\partial}{\partial \theta_j} \log \Pi_{2p} = 2 \frac{\partial}{\partial \theta_j} \Big( \sum_{1 \leq k < \ell \leq p} \log( \cos \theta_k - \cos \theta_\ell) \Big) \\
&&\hspace{-0.5cm}
= - \sum_{j<\ell \leq p} \frac{2 \sin \theta_j}{\cos \theta_j - \cos \theta_\ell} + \sum_{1 \leq k<j} \frac{2 \sin \theta_j}{\cos \theta_k - \cos \theta_j} = \sum_{k \not = j} \frac{2 \sin \theta_j}{\cos \theta_k - \cos \theta_j}, 
\end{eqnarray*}
which leads to \eqref{eq:na_Th_u2p} . 

\smallskip
\noindent
(ii) We remark that, up to an unimportant multiplicative constant, we have
$$ u_{2p+1} (\Theta) = \prod_{1 \leq k<\ell \leq p} \big( \cos \theta_k - \cos \theta_\ell \big)^2 \, \prod_{k=1}^p (1 -\cos \theta_k), $$
and we proceed like in the previous point. 

\smallskip
\noindent
(iii) case $n=2p$: this proof is adapted from that of \cite[Lemma 12.5.2]{faraud2008Analysis}. We have, up to an unimportant multiplicative constant, 
\begin{equation} 
\Pi_{2p}(\Theta) = \prod_{k<\ell} (e^{i \theta_k} + e^{- i \theta_k} - e^{i \theta_\ell} - e^{- i \theta_\ell}).
\label{eq:lem_auxi_prf1}
\end{equation}
We note that the right-hand side of \eqref{eq:lem_auxi_prf1} is the Vandermonde determinant of the unknowns $(e^{i \theta_k} + e^{- i \theta_k})_{k=1}^p$. We recall that if $(X_1, \ldots, X_p)$ are elements of ${\mathbb C}$, the Vandermonde determinant is 
$$ V(X_1, \ldots, X_n) = \textrm{det} \big( (X_k^{p-\ell})_{k\ell} \big) = \prod_{k < \ell} (X_k - X_\ell). $$ 
Thus, we have 
\begin{eqnarray}
\Pi_{2p} (\Theta) &=& V \big( e^{i \theta_1} + e^{- i \theta_1}, \ldots, e^{i \theta_p} + e^{- i \theta_p} \big)  \nonumber\\
&=& \textrm{det} \Big( \big( (e^{i \theta_k} + e^{-i  \theta_k})^{p-\ell} \big)_{k\ell} \Big). \label{eq:Pi2_vandermonde}
\end{eqnarray}
By developing the powers involved in all the terms of the determinant and using elementary manipulations on the rows and columns (see \cite[Eq. (24.39)]{fulton2013representation}), one realizes that 
\begin{eqnarray*}
\Pi_{2p} (\Theta) &=& \textrm{det} \big( ( e^{i (p-\ell) \theta_k} + e^{-i (p-\ell) \theta_k})_{k\ell} \big) \\
&=& \sum_{\sigma \in \mathfrak{S}_p} (-1)^\sigma \prod_{k=1}^p \big( e^{i (p-\sigma(k)) \theta_k} + e^{-i (p-\sigma(k)) \theta_k} \big) , 
\end{eqnarray*}
where $\mathfrak{S}_p$ is the symmetric group of $p$ elements and for $\sigma \in \mathfrak{S}_p$, $(-1)^\sigma$ stands for its signature. Now, taking two derivatives with respect to $\theta_j$, we get 
$$
\frac{\partial^2 \Pi_{2p}}{\partial \theta_j^2} (\Theta) 
= - \sum_{\sigma \in \mathfrak{S}_p} (-1)^\sigma (p-\sigma(j))^2 \prod_{k=1}^p \big( e^{i (p-\sigma(k)) \theta_k} + e^{-i (p-\sigma(k)) \theta_k} \big) . 
$$
Summing over $j$ and realizing that 
$$ \sum_{j=1}^p (p-\sigma(j))^2 = \sum_{j=1}^p (p-j)^2 = \gamma_{2p}, $$
is independent of $\sigma$, we finally get \eqref{eq:Delt_Pin} for $n=2p$. 

\smallskip
\noindent
(iii) case $n=2p+1$: the proof is similar. We note that, up to an unimportant multiplicative constant, we have 
\begin{eqnarray}
\Pi_{2p+1} &=& \prod_{k<\ell} (e^{i \theta_k} + e^{- i \theta_k} - e^{i \theta_\ell} - e^{- i \theta_\ell}) \prod_{k=1}^p (e^{i \theta_k/2} - e^{-i \theta_k/2}) \label{eq:Pi_2p+1_prod}\\
&=&  \textrm{det} \big( ( e^{i (p-\ell + \frac{1}{2}) \theta_k} + e^{-i (p-\ell + \frac{1}{2}) \theta_k})_{k\ell} \big), \nonumber
\end{eqnarray}
where the second inequality again requires elementary manipulations of the determinant (see \cite[Exercise 24.27]{fulton2013representation}. The remainder of the proof is similar. \endproof

We now provide alternate expressions of the radial Laplacian in the following

\begin{proposition} (i) In the case $n=2p$, for any $\varphi \in C^\infty({\mathcal T})$, we have 
\begin{equation}
L \varphi = \sum_{j=1}^p \Big[ \frac{\partial^2 \varphi}{\partial \theta_j^2} + \Big( \sum_{k \not = j} \frac{2}{\cos \theta_k - \cos \theta_j} \Big) \sin \theta_j \frac{\partial \varphi}{\partial \theta_j} \Big]. 
\label{eq:radlap_equiv_2p}
\end{equation}

\smallskip
\noindent
(ii) In the case $n=2p+1$, for any $\varphi \in C^\infty({\mathcal T})$, we have 
\begin{equation}
L \varphi = \sum_{j=1}^p \Big[ \frac{\partial^2 \varphi}{\partial \theta_j^2} + \Big( \sum_{k \not = j} \frac{2}{\cos \theta_k - \cos \theta_j} + \frac{1}{1 - \cos \theta_j} \Big) \sin \theta_j \frac{\partial \varphi}{\partial \theta_j} \Big]. \label{eq:radlap_equiv_2p+1} 
\end{equation}

\smallskip
\noindent
(iii) In both cases $n=2p$ (for $p \geq 2$) and $n=2p+1$ (for $p \geq 1$), we have 
\begin{eqnarray}
L \varphi &=&  \frac{1}{u_n} \nabla_\Theta \cdot \big( u_n \nabla_\Theta \varphi \big), \label{eq:radlap_conserv} \\
&=& \frac{1}{\Pi_n} \Big( \Delta_\Theta + \gamma_n \Big) ( \Pi_n \varphi), \label{eq:radlap_Pi}
\end{eqnarray}
\label{prop:radlap_equiv}
\end{proposition}

\noindent
\textbf{Proof.} (i) and (ii): formulas \eqref{eq:radlap_equiv_2p} and \eqref{eq:radlap_equiv_2p+1} are obtained from \eqref{eq:radlap_2p} and \eqref{eq:radlap_2p+1} respectively, by exchanging $i$ and $j$ inside the second terms of the sums over $i<j$. 

\smallskip
\noindent
(iii) Inserting \eqref{eq:na_Th_u2p} and \eqref{eq:na_Th_u2p+1} into \eqref{eq:radlap_equiv_2p} and \eqref{eq:radlap_equiv_2p+1} respectively, we can write 
$$ L \varphi = \Delta_\Theta \varphi + \frac{1}{u_n} \nabla_\Theta u_n \cdot \nabla_\Theta \varphi = \frac{1}{u_n} \nabla_\Theta \cdot (u_n \nabla_\Theta \varphi), $$
which is \eqref{eq:radlap_conserv}. Equivalently, we can write 
\begin{eqnarray*} 
L \varphi &=& \frac{1}{\Pi_n^2} \nabla_\Theta \cdot (\Pi_n^2 \nabla_\Theta \varphi) = \frac{1}{\Pi_n} \Big( \nabla_\Theta \Pi_n \cdot \nabla_\Theta \varphi + \nabla_\Theta \cdot (\Pi_n \cdot \nabla_\Theta \varphi) \Big)\\
&=& \frac{1}{\Pi_n} \Big( \nabla_\Theta \Pi_n \cdot \nabla_\Theta \varphi + \nabla_\Theta \cdot \big( \nabla_\Theta (\Pi_n \varphi) - \varphi \nabla_\Theta \Pi_n \big) \Big) \\
&=& \frac{1}{\Pi_n} \Big( \Delta_\Theta (\Pi_n \varphi) - \varphi \Delta_\Theta (\Pi_n) \Big), 
\end{eqnarray*}
which, with \eqref{eq:Delt_Pin}, leads to \eqref{eq:radlap_Pi}. \endproof

\subsection{Spectral theory and inversion of the radial Laplacian}
\label{subsec_spectral}

In this and the following sections, we denote by $G$ either of the groups $\textrm{SO}_{2p}{\mathbb R}$ or $\textrm{SO}_{2p+1}{\mathbb R}$. We state a few facts concerning the spectral decomposition of the Laplacian $-\Delta$ and the radial Laplacian $-L$. 

We first start with the Laplacian $- \Delta$. We know \cite{gallot1990riemannian} that the spectrum of $- \Delta$ on a Riemannian manifold consists of an increasing sequence of discrete eigenvalues $\lambda_1 < \lambda_2 < \ldots \lambda_n \ldots$, tending to $+\infty$, that the associated eigenspaces are finite dimensional and that their direct sum is dense in $L^2(G)$. Furthermore, $\lambda_1 = 0$ and its eigenspace is one-dimensional and consists of the constant functions. However, in the case of the Laplacian on $G$, we can be more specific about the eigenspaces. We first need to introduce a few concepts from representation theory. 

Here, we refer to \cite[Ch. 6]{faraud2008Analysis}. A representation $(V,\pi)$ of $G$ is a pair with $V$ a vector space over ${\mathbb C}$ and $\pi$ a group morphism $G \to \textrm{Aut}(V)$ where $\textrm{Aut}(V)$ is the group of linear automorphisms of $V$. If $W$ is a subspace of $V$ such that $\pi(g) W =W$, $\forall g \in G$, then $(W, \tilde \pi)$ where $\tilde \pi(g) = \pi(g)|_{W}$ is a representation called a subrepresentation of $V$. An irreducible representation is a representation which has no proper subrepresentation. When $G$ is compact, any irreducible representation is finite dimensional and any finite-dimensional representation can be decomposed in the direct sum of irreducible representations. Two representations $(V,\pi)$ and $(V', \pi')$ are equivalent if and only if there exists a linear isomorphism $A$: $V \to V'$ which intertwines the two representations, i.e. $A \pi(g) = \pi'(g) A$, $\forall g \in G$. We denote by $\hat G$ the set of equivalence classes of irreducible representations of $G$. For any $\lambda \in \hat G$, let us choose a representative $(V_\lambda, \pi_\lambda)$ in $\lambda$. We denote by ${\mathcal M}_\lambda$ the space of functions $f$: $G \to {\mathbb C}$ such that there exists $A \in \textrm{Aut}(V_\lambda)$ and $f(g) = \textrm{Tr} \{ A \pi_\lambda(g) \}$. Then, the Peter-Weyl theorem asserts that 
\begin{equation} 
L^2(G) = \overline{ \bigoplus_{\lambda \in \hat G} {\mathcal M}_\lambda}, 
\label{eq:Peter-Weyl}
\end{equation}
where $L^2(G)$ is the space of square integrable functions on $G$. Furthermore \cite[Prop. 8.2.1]{faraud2008Analysis}, all functions $f \in {\mathcal M}_\lambda$ are eigenfunctions of the Laplace operator associated with the same eigenvalue. So, the decomposition \eqref{eq:Peter-Weyl} is the eigenspace decomposition of the Laplace operator on $G$. 

Thanks to this, we can easily solve the following problem: given $h \in L^2(G)$, find $f$ such that 
\begin{equation}
- \Delta f = h, 
\label{eq:Delf=h}
\end{equation}
Obviously, $h$ must be chosen to satisfy the compatibility equation 
\begin{equation} 
\int_G h(x) \, dx = 0, 
\label{eq:compat}
\end{equation}
where we now denote the normalized Haar measure on $G$ by $dx$. A usual method is to use a variational formulation. Supposing $f \in H^2(G)$ (where $H^k(G)$ is defined recursively by $H^k(G) = \{ f \in H^{k-1}(G) \, \, | \, \, \nabla f \in H^{k-1}(G) \}$, and $H^0(G) = L^2(G)$) and taking a test function $v \in H^1(G)$, multiplying \eqref{eq:Delf=h} by $v$ and using Stokes formula, we get 
\begin{equation}
\int_G \nabla f \cdot \nabla v \, dx = \int_G h \, v \, dx, \quad \forall v \in H^1(G). 
\label{eq:var_form_1}
\end{equation}
Thanks to \eqref{eq:compat}, Eq. \eqref{eq:var_form_1} is equivalent to the same formulation but with $v$ taken in $H^1_0(G)$ where $H^1_0(G) = \{ v \in H^1(G) \, \, | \, \, \int_G v(x) \, dx = 0 \}$. Now, since clearly, a solution $f$ to \eqref{eq:Delf=h} is defined up to constants, it is reasonable to search for a unique $f$ in the space $H^1_0(G)$. Indeed, we have the following 

\begin{theorem}
Given $h \in L^2(G)$ satisfying the compatibility condition \eqref{eq:compat}, the variational formulation 
\begin{eqnarray}
&&\hspace{-1cm}
f \in H^1_0(G), \label{eq:var_form_2-1} \\
&&\hspace{-1cm}
\int_G \nabla f \cdot \nabla v \, dx = \int_G h \, v \, dx, \quad \forall v \in H^1_0(G), \label{eq:var_form_2-2}
\end{eqnarray}
has a unique solution, referred to as the variational solution of the problem $- \Delta f = h$. Furthermore, the $H^1$ norm of $f$ is controlled by the $L^2$ norm of $h$. 
\label{thm_inversion_lap}
\end{theorem}

\noindent
\textbf{Proof.} By the Lax-Milgram theorem, it is enough to show a Poincare inequality. But because there is a spectral gap, i.e. $\lambda_2 >0$, we have $$ \int_G |\nabla f|^2 \, dx \geq \lambda_2 \int_G |f|^2 \, dx, \quad \forall f \in H^1_0(G). $$
This shows the result. \endproof

\begin{remark} 
Thanks to the Peter-Weyl theorem and a Plancherel theorem \cite[Thm 6.4.2]{faraud2008Analysis}, we could express the Fourier coefficients of $f$ in terms of the Fourier coefficients of $h$ and of the eigenvalues of $\Delta$, where the Fourier coefficient $\hat f(\lambda)$ is defined as an endomorphism of $V_\lambda$ given by $\hat f(\lambda) = \int_G f(x) \, \pi_\lambda(x^{-1}) \, dx$. \end{remark}

Now, we focus on class functions. The character of a representation $(V,\pi)$ is the map $G \to {\mathbb C}$ defined by $\chi_\pi(g) = \textrm{Tr}\{ \pi(g) \}$, $\forall g \in G$. It is a class function. Now, let $\lambda \in \hat G$ and denote by $\chi_\lambda$ the character of any representative $(V_\lambda, \pi_\lambda)$ of $\lambda$. Then, the subspace of ${\mathcal M}_\lambda$ consisting of class functions is one-dimensional and spanned by $\chi_\lambda$. The system $\{ \chi_\lambda \}_{\lambda \in \hat G}$ is a Hilbert basis of $L^2(G)_{\textrm{class}}$, the subspace of $L^2(G)$ consisting of class functions \cite[Prop. 6.5.3]{faraud2008Analysis} (the fact that $\{ \chi_\lambda \}_{\lambda \in \hat G}$ is an orthonormal system is known as Schur's orthogonality relations). Now, since $\chi_\lambda$ is an eigenfunction of $\Delta$ and is a class function, $\varphi_\lambda =: \varphi_{\chi_\lambda}$ is an eigenfunction of the radial Laplacian $L$ (we recall that for a class function $f$, $\varphi_f$ is defined by \eqref{eq:def_phif}). 

We introduce $L^2({\mathcal T})$, the set of functions $\varphi$: ${\mathcal T} \to {\mathbb R}$, such that $\int_{{\mathcal T}} |\varphi(\Theta)|^2 \, u_n(\Theta) \, d \Theta < \infty$ and $H^k({\mathcal T})$, defined recursively by $\varphi \in H^k({\mathcal T}) \Longleftrightarrow (\varphi \in H^{k-1}({\mathcal T})$ and $\nabla_\Theta \varphi \in H^{k-1}({\mathcal T}))$, with $H^0({\mathcal T}) = L^2({\mathcal T})$. With these notations, $(\varphi_\lambda)_{\lambda \in \hat G}$ is a Hilbert basis of $L^2({\mathcal T})$ consisting of eigenfunctions of $L$. It follows that $L$ has the same eigenvalues as $\Delta$. Consequently, analogously to the case of the Laplacian, we can define a variational solution of the problem $-L \varphi = \eta$. Let us define $H^1_0({\mathcal T}) = \{ v \in H^1({\mathcal T}) \, \, | \, \, \int_{{\mathcal T}} v(\Theta) \, u_n(\Theta) \, d \Theta = 0 \}$. Then, we have 

\begin{theorem}
Given $\eta \in L^2({\mathcal T})$ satisfying the compatibility condition
\begin{equation}
\int_{{\mathcal T}} \eta(\Theta) \, u_n(\Theta) \, d \Theta = 0, 
\label{eq:compat_L}
\end{equation}
the variational formulation 
\begin{eqnarray}
&&\hspace{-1cm}
\varphi \in H^1_0({\mathcal T}), \label{eq:var_form_L-1} \\
&&\hspace{-1cm}
\int_{{\mathcal T}} \nabla_\Theta \varphi \cdot \nabla_\Theta v \, u_n \, d\Theta = \int_{{\mathcal T}} \eta \, v \, u_n \,  d\Theta, \quad \forall v \in H^1_0({\mathcal T}), \label{eq:var_form_L-2}
\end{eqnarray}
has a unique solution referred to as the variational solution of the problem $-L \varphi = \eta$. Furthermore, the $H^1$ norm of $\varphi$ is controlled by the $L^2$ norm of $\eta$. 
\label{thm_inversion_radlap}
\end{theorem}

Similarly to the case of the inversion of $-\Delta$, we can also recover the Fourier coefficients of $\varphi$ (in the Hilbert basis $\{\varphi_\lambda\}_{\lambda \in \hat G}$) from those of $\eta$ and the eigenvalues of the Laplace operator. The next section is devoted to the computation of the eigenvalues of $\Delta$ and $L$.

\subsection{Characters and eigenvalues of the Laplacian}
\label{subsec_characters}

As seen in the previous section, characters are at the same time eigenfunctions of the Laplacian and  class functions. Expressions of the characters on the maximal torus ${\mathbb T}$ are available thanks to the Weyl character formula. Applying the radial Laplacian to them, we will deduce expressions of the eigenvalues of the Laplacian. All this is expressed in the following theorem: 

\begin{theorem} (i) Case $n=2p$ for $p \geq 2$. The irreducible representations of $\textrm{SO}_{2p}{\mathbb R}$ are indexed by a $p$-tuple of integers $\lambda = (\lambda_1, \ldots, \lambda_p) \in {\mathbb Z}^p$ such that 
\begin{equation}
\lambda_1 \geq \lambda_2 \geq \ldots \geq \lambda_{p-1} \geq |\lambda_p| \geq 0.
\label{eq:dominant_2p}
\end{equation}
 Let $\chi_\lambda$ be its character and define $\varphi_\lambda =: \varphi_{\chi_\lambda}$ with $\varphi_f$ defined by \eqref{eq:def_phif}. Then $\varphi_\lambda$ is given by the Weyl character formula \cite[Eq. (24.40)]{fulton2013representation}: 
\begin{equation}
\varphi_\lambda (\Theta) = \frac{ \textrm{det} \, \big( (e^{i \theta_k (\lambda_\ell + \delta_\ell)} + e^{-i \theta_k (\lambda_\ell + \delta_\ell)})_{k \ell} \big)}{\Pi_{2p}(\Theta)} + \frac{ \textrm{det} \, \big( (e^{i \theta_k (\lambda_\ell + \delta_\ell)} - e^{-i \theta_k (\lambda_\ell + \delta_\ell)})_{k \ell} \big)}{\Pi_{2p}(\Theta)},
\label{eq:chi_lambda_2p}
\end{equation}
with $ \delta_\ell = p-\ell$, $\forall \ell \in \{1, \ldots, p\}$. We denote by $\delta$ the $p$-tuple $\delta = (\delta_\ell)_{\ell=1}^p$.

\smallskip
\noindent
(ii) Case $n=2p+1$ for $p \geq 1$. The irreducible representations of $\textrm{SO}_{2p+1}{\mathbb R}$ are indexed by a $p$-tuple of integers $\lambda = (\lambda_1, \ldots, \lambda_p) \in {\mathbb Z}^p$ such that 
\begin{equation}
\lambda_1 \geq \lambda_2 \geq \ldots \geq \lambda_{p-1} \geq \lambda_p \geq 0. 
\label{eq:dominant_2p+1}
\end{equation}
Then, the Weyl character formula \cite[Eq. (24.28)]{fulton2013representation} is written in this case: 
\begin{equation}
\varphi_\lambda (\Theta) = \frac{ \textrm{det} \, \big( (e^{i \theta_k (\lambda_\ell + \delta_\ell)} - e^{-i \theta_k (\lambda_\ell + \delta_\ell)})_{k \ell} \big)}{\Pi_{2p+1}(\Theta)},
\label{eq:chi_lambda_2p+1}
\end{equation}
with $\delta_\ell = p-\ell + \frac{1}{2}$, $\forall \ell \in \{1, \ldots, p\}$. We denote by $\delta$ the $p$-tuple $\delta = (\delta_\ell)_{\ell=1}^p$. 

\smallskip
\noindent
(iii) In both the cases  $n=2p$ for $p \geq 2$ and $n=2p+1$ for $p \geq 1$, $\varphi_\lambda$ and $\chi_\lambda$ are eigenfunctions of $-L$ and $- \Delta$ respectively, associated with the eigenvalue $\kappa_\lambda$ given by
\begin{equation}
\kappa_\lambda = \| \lambda + \delta \|^2 - \| \delta \|^2, 
\label{eq:kappa_lambda}
\end{equation}
where $\| \cdot \|$ is the Euclidean norm in ${\mathbb R}^p$. Note that $\| \delta \|^2 = \gamma_{2p}$ or $\gamma_{2p+1}$ according to the cases considered (see Eq. \eqref{eq:gamma}). 

\label{thm:characters}
\end{theorem}

\noindent
\textbf{Proof.} (i) \& (ii) lie at the heart of Lie group representation theory and we refer the reader to e.g. \cite{fulton2013representation}. Here, we just make sure that Eqs. \eqref{eq:chi_lambda_2p} and \eqref{eq:chi_lambda_2p+1} make sense, i.e. that they define functions $\varphi_\lambda$ in ${\mathcal T}$, in other words that $\varphi_\lambda$ are $2\pi$ periodic in each variable~$\theta_j$ and invariant by the Weyl group. 

\smallskip
\noindent
$2 \pi$-periodicity follows immediately from the assumption that $\lambda \in {\mathbb Z}^p$. Indeed, in the case $n=2p$, changing one of the $\theta_j$ into $\theta_j + 2 \pi$ leaves both the numerator and the denominator of \eqref{eq:chi_lambda_2p} invariant. In the case $n=2p+1$, both change sign so that their quotient \eqref{eq:chi_lambda_2p+1} remains invariant. In fact, it can be shown that $\lambda \in {\mathbb Z}^p$ is also a necessary condition for $\varphi_\lambda$ to be $2 \pi$-periodic. 

\smallskip
\noindent
Invariance by the Weyl group will be shown by showing the following stronger property, which is an interesting observation in its own right:  

\begin{itemize}
\item Case $n=2p+1$: there exists a symmetric polynomial $P(X_1, \ldots, X_p)$ such that 
\begin{equation} 
\varphi_\lambda(\Theta) = P (e^{i \theta_1} + e^{-i \theta_1}, \ldots, e^{i \theta_p} + e^{-i \theta_p}), 
\label{eq:sym_pol_2p+1}
\end{equation}

\item Case $n=2p$: there exist two symmetric polynomials $P(X_1, \ldots, X_p)$ and $Q(X_1, \ldots, X_p)$ such that 
\begin{equation} 
\hspace{-0.5cm}
\varphi_\lambda(\Theta) = P (e^{i \theta_1} + e^{-i \theta_1}, \ldots, e^{i \theta_p} + e^{-i \theta_p}) + \prod_{k=1}^n (e^{i \theta_k} - e^{-i \theta_j}) Q (e^{i \theta_1} + e^{-i \theta_1}, \ldots, e^{i \theta_p} + e^{-i \theta_p}). 
\label{eq:sym_pol_2p}
\end{equation}
\end{itemize}

It is clear that \eqref{eq:sym_pol_2p} and \eqref{eq:sym_pol_2p+1} imply invariance by the Weyl group (see discussion after \eqref{eq:def_phif} for the description of the Weyl group for even and odd orthogonal groups). 

In the case $n=2p+1$, we first focus on the numerator of  \eqref{eq:chi_lambda_2p+1}. We have 
\begin{eqnarray*} 
e^{i \theta_k (\lambda_\ell + \delta_\ell)} - e^{-i \theta_k (\lambda_\ell + \delta_\ell)} &=& e^{i \theta_k (\lambda_\ell + p- \ell + \frac{1}{2})} - e^{-i \theta_k (\lambda_\ell + p- \ell + \frac{1}{2})} \\
&=& (e^{i \frac{\theta_k}{2}} - e^{-i \frac{\theta_k}{2}}) \Big( 1 + \sum_{m=1}^{\lambda_\ell + p-\ell} (e^{i m \theta_k} + e^{-i m \theta_k}) \Big). 
\end{eqnarray*}
It is easy to show by induction that there exists a polynomial $p_m$ of degree $m$ such that 
$$ e^{i m \theta_k} + e^{-i m \theta_k} = p_m (e^{i \theta_k} + e^{-i \theta_k}). $$
Thus, there exists a polynomial $q_m$ of degree $m$ such that 
$$ e^{i \theta_k (\lambda_\ell + \delta_\ell)} - e^{-i \theta_k (\lambda_\ell + \delta_\ell)} = 
 (e^{i \frac{\theta_k}{2}} - e^{-i \frac{\theta_k}{2}}) \, q_{\lambda_\ell + p-\ell}(e^{i \theta_k} + e^{-i \theta_k}). $$
It follows that 
\begin{eqnarray} 
&&\hspace{-1cm}
\textrm{det} \, \big( (e^{i \theta_k (\lambda_\ell + \delta_\ell)} - e^{-i \theta_k (\lambda_\ell + \delta_\ell)})_{k \ell} \big) \nonumber \\
&&\hspace{3cm}
= \Big( \prod_{k=1}^p (e^{i \frac{\theta_k}{2}} - e^{-i \frac{\theta_k}{2}}) \Big) \, \textrm{det} \, \Big( \big(q_{\lambda_\ell + p-\ell}(e^{i \theta_k} + e^{-i \theta_k}) \big)_{k \ell} \Big) \label{eq:num_phi_lam_2p+1}\\
&&\hspace{3cm}
= \Big( \prod_{k=1}^p (e^{i \frac{\theta_k}{2}} - e^{-i \frac{\theta_k}{2}}) \Big) R(e^{i \theta_1} + e^{-i \theta_1}, \ldots, e^{i \theta_p} + e^{-i \theta_p}), \nonumber
\end{eqnarray}
where $R(X_1, \ldots, X_p)$ is an alternating polynomial (i.e. it is changed in its opposite when $X_j$ and $X_k$ are exchanged, for any pair $j \not = k$). Thanks to \eqref{eq:Pi_2p+1_prod}, we get 
\begin{equation}
\varphi_\lambda(\Theta) = \frac{R(e^{i \theta_1} + e^{-i \theta_1}, \ldots, e^{i \theta_p} + e^{-i \theta_p})}{\prod_{k<\ell} (e^{i \theta_k} + e^{- i \theta_k} - e^{i \theta_\ell} - e^{- i \theta_\ell})}. \label{eq:phi_lam_poly_2p+1}
\end{equation}
From \eqref{eq:num_phi_lam_2p+1}, it is clear that any factor $(X_j-X_k)$ for $j<k$ divides $R$. Since all these factors are mutually prime polynomials, the denominator divides the numerator. Hence, the right-hand side of \eqref{eq:phi_lam_poly_2p+1} is a polynomial. Since both the numerator and denominator of \eqref{eq:phi_lam_poly_2p+1} are alternating, their quotient is symmetric, which ends the proof in the case $n=2p+1$. 

The case $n=2p$ is similar. The first term at the right-hand side of \eqref{eq:chi_lambda_2p} gives rise to the first term in \eqref{eq:sym_pol_2p}. Indeed, the numerator is clearly an alternating polynomial of $(e^{i \theta_k} + e^{-i \theta_k})_{k=1}^p$ and the denominator too thanks to \eqref{eq:Pi2_vandermonde}. Furthermore, each factor of the denominator divides the numerator and are mutually prime, so the quotient is a symmetric polynomial. The second term at the right-hand side of \eqref{eq:chi_lambda_2p} gives rise to the second term in \eqref{eq:sym_pol_2p} thanks to similar computations as those made in the case $n=2p+1$. 

\smallskip
\noindent

\smallskip
\noindent
(iii) The general methodology is outlined in \cite[\S 5.7]{gurarie2007symmetries} but the proof is quite abstract. On the other hand, proving \eqref{eq:kappa_lambda} is quite easy, using \eqref{eq:radlap_Pi}. Let us consider the case $n=2p$ for instance. With \eqref{eq:chi_lambda_2p}, we have 
$$ \Pi_{2p} \varphi_\lambda (\Theta) = \sum_{\sigma \in \mathfrak{S}_p} (-1)^\sigma \prod_{j=1}^p \Big( e^{i \ell_{\sigma(j)} \theta_j} + e^{-i \ell_{\sigma(j)} \theta_j} \Big) + \sum_{\sigma \in \mathfrak{S}_p} (-1)^\sigma \prod_{j=1}^p \Big( e^{i \ell_{\sigma(j)} \theta_j} - e^{-i \ell_{\sigma(j)} \theta_j} \Big). $$
Now, like in the proof of Lemma \ref{lem_auxi_un_Pin} (iii), we have 
$$ - \Delta_\Theta ( \Pi_{2p} \varphi_\lambda) = \Big( \sum_{k=1}^p \ell_k^2 \Big) \Pi_{2p} \varphi_\lambda. $$
Then, with \eqref{eq:radlap_Pi} and \eqref{eq:Delt_Pin}, we get \eqref{eq:kappa_lambda}. The case $n=2p+1$ is similar. \endproof

We end this section with some terminology. The $p$-tuple $\lambda$ lies on the lattice ${\mathbb Z}^p$, which is called the weight lattice. Any $\lambda \in {\mathbb Z}^p$ is called a weight. The weights satisfying the restrictions \eqref{eq:dominant_2p} or \eqref{eq:dominant_2p+1} are called dominant weights, and the regions of ${\mathbb R}^{p}$ defined by the inequalities \eqref{eq:dominant_2p} or \eqref{eq:dominant_2p+1} are called the closed Weyl chambers. All irreducible representations of $\textrm{SO}_n{\mathbb R}$ induce irreducible representations of the Lie algebra $\mathfrak{so}_n{\mathbb R}$. The converse is not true. The weight lattice of $\mathfrak{so}_n{\mathbb R}$ is equal to ${\mathbb Z}^p \cup (\frac{1}{2} + {\mathbb Z})^p$. The Weyl chambers of $\mathfrak{so}_n{\mathbb R}$ and $\textrm{SO}_n{\mathbb R}$ on the other hand are the same. All representations of $\mathfrak{so}_n{\mathbb R}$ lift to representations of the spin group $\textrm{Spin}_n{\mathbb R}$ which is a simply connected doubly-sheeted covering of $\textrm{SO}_n{\mathbb R}$ and is such that $\textrm{SO}_n{\mathbb R} = \textrm{Spin}_n{\mathbb R} / \{\pm1\}$. Half-integer weights correspond to representations of $\textrm{Spin}_n{\mathbb R}$ that are not trivial on $\{\pm1\}$ and thus, cannot induce representations of $\textrm{SO}_n{\mathbb R}$ (but some authors call these representations doubled-valued representations of $\textrm{SO}_n{\mathbb R}$).

\setcounter{equation}{0}
\section{Conclusion}
\label{sec_conclu}

In this work, we have reviewed elements of the spectral theory of the Laplacian on rotation groups. One key element is to find the expression of the radial Laplacian, which is the restriction of the Laplacian on the maximal torus. Adequate expressions of this radial Laplacian, together with the use of the Weyl chararcter formula lead to explicit expressions of the eigenvalues of the Laplacian. Although the material presented in these notes is familiar to Lie group experts, it is not easily found in the literature in a concise and accessible way. We hope these notes will prove useful to those working with rotation groups without being experts of Lie group theory.

\appendix

\setcounter{equation}{0}
\section{Proof of Lemma \ref{lem:vol_form}}
\label{secapp_proportional}

By uniqueness of the Haar measure, to show that there is a positive constant $C$ such that $d \omega = C d\mu$, it is enough to show that the Riemannian volume form is left invariant. First, we define an orientation of $\mathfrak{so}_n{\mathbb R}$ associated with the lexicographic ordering of the orthonormal basis $(F_{ij})_{1 \leq i < j \leq n}$. For $A \in \textrm{SO}_n{\mathbb R}$, the tangent space $T_A$ is oriented by the same lexicographic ordering of the basis $(AF_{ij})_{i < j}$, and this gives rise to a global orientation of $\textrm{SO}_n{\mathbb R}$. We note that $(AF_{ij})_{i < j}$ is an orthonormal basis of $T_A$ precisely because $A$ is an orthogonal matrix. Then, we construct the Riemannian volume form of $\textrm{SO}_n{\mathbb R}$ associated with this orientation, following the construction of \cite[Section 4.10]{warner1983foundations}. We denote by $(\omega_{ij}(A))_{i < j}$ the dual basis of $(AF_{ij})_{i < j}$. The inner-product on $T_A$ generates a pairing between $T_A$ and its dual $T_A^*$. This pairing is such that the image $\omega_X \in T_A^*$ of an element $X \in T_A$ is such that $\langle \omega_X, Y \rangle = X \cdot Y$ for all $Y \in T_A$, where the bracket indicates duality between $T_A^*$ and $T_A$. In turn, this pairing induces an inner product on $T_A^*$, such that $\omega_X \cdot \omega_Y = X \cdot Y$. Thus, since $(AF_{ij})_{i < j}$ is an orthonormal basis of $T_A$, $(\omega_{ij}(A))_{i < j}$ is an orthonormal basis of $T_A^*$, so $\omega (A) = \bigwedge_{i<j} \omega_{ij}(A)$ is the volume form where $\bigwedge_{i<j}$ stands for the exterior product of the forms $\omega_{ij}$ in lexicographic order. 

Now, we show that this volume form is left invariant. Let $\sigma \in \textrm{SO}_n{\mathbb R}$ and $\ell_\sigma$ the morphism of $\textrm{SO}_n{\mathbb R}$ such that $\ell_\sigma (A) = \sigma A$. The form $\omega$ is left invariant if and only if $\delta \ell_\sigma (\omega) = \omega$ for all $\sigma \in \textrm{SO}_n{\mathbb R}$, where $\delta \ell_\sigma (\omega)$ is the pullback of $\omega$ by $\ell_\sigma$. Let $(AX_{ij})_{i<j}$, with $X_{ij} \in \mathfrak{so}_n{\mathbb R}$ for all $i<j$ be elements of $T_A$. Then, $\delta \ell_\sigma (\omega)$ is defined by 
\begin{eqnarray*} 
&&\hspace{-1cm}
\delta \ell_\sigma (\omega) (A) \big( AX_{12}, AX_{13}, \ldots, AX_{n-1 \, n} \big) \\ 
&&\hspace{1cm}
= \omega(\sigma A) \big( d\ell_\sigma(A) (AX_{12}), d\ell_\sigma(A) (AX_{13}), \ldots, d\ell_\sigma(A) (AX_{n-1 \, n}) \big), 
\end{eqnarray*}
where $d\ell_\sigma(A)$ is the differential of $\ell_\sigma$ at $A$. It is readily checked that $d\ell_\sigma(A) (AX) = \sigma A X$, for any $X \in \mathfrak{so}_n{\mathbb R}$. Thus, 
\begin{eqnarray*} 
\delta \ell_\sigma (\omega) (A) \big( AX_{12}, AX_{13}, \ldots, AX_{n-1 \, n} \big) &=& \omega(\sigma A) (\sigma A X_{12}, \sigma A X_{13}, \ldots, \sigma A X_{n-1 \, n} ) \\
&=& \sum_{\tau \in \mathfrak{S}_N} \varepsilon(\tau) \prod_{i<j} \omega_{ij}(\sigma A) \big( \sigma A X_{\tau(i,j)} \big), 
\end{eqnarray*}
where $N = \frac{n(n-1)}{2}$ is the dimension of $T_A$, $\mathfrak{S}_N$ is the permutation group of $N$ elements, $\varepsilon(\tau)$ is the signature of the permutation $\tau$ and where the elements of $\{1, \ldots, N\}$ are labeled by the pairs $(i,j)$ with $1 \leq i < j \leq n$. Hence, $\tau(i,j)$ denotes the image of the element $(i,j)$ of $\{1, \ldots, N\}$ by the permutation $\tau$. 

Now, if we show that 
\begin{equation}
\omega_{ij}(\sigma A) \big( \sigma A X \big) = \omega_{ij}(A) \big(A X \big), \, \,  \forall A, \, \sigma \in \textrm{SO}_n{\mathbb R}, \, \,  \forall X \in \mathfrak{so}_n{\mathbb R}, \, \,  \forall (i,j) \, | \,  1 \leq i < j \leq n,  
\label{eq:omij}
\end{equation}
then, we obviously deduce that $\delta \ell_\sigma(\omega) = \omega$, which concludes the proof. To prove \eqref{eq:omij}, it is enough to prove it for all basis vectors of a basis of $\mathfrak{so}_n{\mathbb R}$. We take the basis $(F_{k \ell})_{k < \ell}$. We have $ \omega_{ij}(\sigma A) \big( \sigma A F_{k \ell} \big) = \delta_{ij} \delta_{k \ell}$ because $(\omega_{ij}(\sigma A))_{i<j}$ is the dual basis of $(\sigma A F_{k \ell})_{k < \ell}$. But this being true for all $\sigma$, it is true for $\sigma = \textrm{I}$ which gives $\omega_{ij}(A) \big(A F_{k \ell} \big) = \delta_{ij} \delta_{k \ell}$. Hence, $\omega_{ij}(\sigma A) \big( \sigma A F_{k \ell} \big) = \omega_{ij}(A) \big(A F_{k \ell} \big)$, which shows \eqref{eq:omij}. \endproof


\setcounter{equation}{0}
\section{Proof of Lemma \ref{lem:laplacians=}}
\label{secapp_lem:laplacians=_proof}

From \eqref{eq:PTA_form} and \eqref{eq:rhoX_prop}, we have, for any $C^\infty$ function $f$ on $\textrm{So}_n{\mathbb R}$ and any $A \in \textrm{So}_n{\mathbb R}$: 
$$ (P_{T_A} E_{ij}) (f) (A)  = \Big( \rho \big( \frac{A^T E_{ij} - E_{ij}^T A}{2} \big) f \Big) (A).
$$
A straightforward computation shows that 
$$ A^T E_{ij} - E_{ij}^T A =  \sum_{k=1}^n A_{i k} F_{j k}, $$
where $A_{ik}$ denote the entries of the matrix $A$. Since $\rho(X)$ is linear with respect to $X$, we have 
$$ \rho ( A^T E_{ij} - E_{ij}^T A ) = \sum_{k=1}^n A_{i k} \rho \big(F_{jk} \big), $$
and 
$$ \big( \rho ( A^T E_{ij} - E_{ij}^T A ) \big)^2 = \sum_{k, \ell} \, A_{i k} A_{i m} \, \rho \big(F_{jk} \big) \, \rho \big(F_{jm} \big). $$
Now, with \eqref{eq:lapriem_def}, we get 
$$ \Delta_M f = \frac{1}{2} \sum_{i, j, k, \ell } A_{i k} A_{i m} \, \rho \big(F_{jk} \big) \, \rho \big(F_{jm} \big) = \frac{1}{2} \sum_{j,k} \rho \big(F_{jk} \big)^2 = \sum_{1 \leq j < k \leq n} \rho \big(F_{jk} \big)^2 = \Delta_G f, $$
where, in the second equality, we have used that $A$ is orthogonal,  in third equality, that $F_{ji} = - F_{ih}$ and in the fourth one, Eq. \eqref{eq:laplie_def}. This ends the proof. \endproof


\end{document}